\documentclass[preprint,12pt,numbers,numbers,sort&compress]{elsarticle}

\usepackage{graphicx}
\usepackage{dcolumn}
\usepackage{bm}
\usepackage[utf8]{inputenc}
\usepackage[T1]{fontenc}
\usepackage{subcaption}
\usepackage{graphicx}

\makeatletter
\def\@email#1#2{%
 \endgroup
 \patchcmd{\titleblock@produce}
  {\frontmatter@RRAPformat}
  {\frontmatter@RRAPformat{\produce@RRAP{*#1\href{mailto:#2}{#2}}}\frontmatter@RRAPformat}
  {}{}
}%
\makeatother
\usepackage{amssymb}                                            
\usepackage{mathtools}
\usepackage{physics}
\usepackage{csquotes}   
\MakeOuterQuote{"}

\numberwithin{equation}{section}
\numberwithin{figure}{section}

\newcommand{\rmi}{\mathrm{i}}
\def\defeq{\mathrel{\mathop:}=}
\newcommand{\DGamma}{\mathsf{D}_{\Gamma}}

\newcommand{\Reals}{\mathbb{R}}
\newcommand{\Complex}{\mathbb{C}}
\newcommand{\JJ}{\mathbb{J}}
\renewcommand{\SS}{\mathbb{S}}
\DeclareMathOperator{\Disc}{Disc}
\newcommand{\vone}{\vb{1}}
\newcommand{\ufrak}{\mathfrak{u}}

\newcommand\SO{\mathsf{SO}}
\newcommand\SE{\mathsf{SE}}

\newcommand\parentheses[1]{\!\left(#1\right)}
\newcommand\brackets[1]{\!\left[#1\right]}

\DeclarePairedDelimiterX\pair[2]{\langle}{\rangle}{#1,#2}
\DeclarePairedDelimiterX\poisson[2]{\{}{\}}{#1,#2}
\newcommand{\Tee}{\mathsf{T}}

\newcommand{\J}{\mathbf{J}}
\newcommand{\z}{\mathbf{z}}
\newcommand{\Z}{\mathbf{Z}}
\renewcommand{\r}{\mathbf{r}}

\DeclareMathOperator{\sign}{sign}
\DeclareMathOperator{\diag}{diag}
\newcommand{\Equilibrium}{\mathcal{E}}
\newcommand{\Singularity}{\mathcal{S}}
\newcommand{\EqTriPM}{\Equilibrium_{\mathrm{tri}}^\pm}

\newcommand{\SingTriple}{\Singularity_{\mathrm{triple}}}
\newcommand{\res}[3]{\mathrm{res}\left(#1,#2;#3\right)}

\newcommand{\sfx}{\mathsf{x}}
\newcommand{\sfy}{\mathsf{y}}
\newcommand{\sfP}{\mathsf{P}}

\journal{Physica D}

\begin{document}

\begin{frontmatter}

\title{Phase portraits and the bifurcation set for the three-vortex interaction system}

\author{Atul Anurag} 
\affiliation{organization={Department of Mathematics, Ramapo College of New Jersey},%
            addressline={505 Ramapo Valley Rd.}, 
            city={Mahwah},
            postcode={07430}, 
            state={NJ},
            country={USA}}

\author{Roy H. Goodman} %

\affiliation{organization={Department of Mathematical Sciences, New Jersey Institute of Technology},%
            addressline={323 Martin Luther King Blvd.}, 
            city={Newark},
            postcode={NJ}, 
            state={07102},
            country={USA}}

\begin{abstract}
We derive a symplectic reduction of the evolution equations for a system of three interacting point vortices in the plane, first introducing Jacobi coordinates, then Lie-Poisson reductions, and finally reparameterizing the resulting leaves, arriving at an integrable system on a topologically nontrivial phase space surface. The reduced system is convenient for describing all aspects of three-vortex dynamics, including finite-time collapse, the calculation of relative equilibria and their stability, and scattering. We use the final simplified system to succinctly and geometrically explain a kind of bifurcation diagram that has appeared in the literature.
\end{abstract}

\begin{highlights}
\item We describe a multi-stage Hamiltonian reduction of the three-vortex interaction, first introducing Jacobi coordinates, second, a Lie-Poisson reduction, and finally a reparameterization of the symplectic leaves, arriving at an integrable system on a topologically non-trivial phase space surface.
\item The reduced system is convenient for describing all aspects of three vortex dynamics, including finite-time collapse, the enumeration of relative equilibria and their stability, and scattering.
\item We describe a bifurcation diagram that explains the global phase portraits as the vortex circulations vary.
\end{highlights}

\end{frontmatter}

\section{Introduction}

This paper introduces a Hamiltonian reduction of the system governing the motion of three interacting point vortices in the plane. This builds on our recent approach in Reference~\cite{Anurag:2024}, eliminating the need to consider two cases separately. The modern machinery used here yields a reduced system in a straightforward form, free of singularities introduced by earlier reductions. The resulting equations prove helpful in understanding all aspects of the dynamics. They are based on a generalized form of Hamiltonian reduction by stages: at each step, the equations maintain their Hamiltonian structure. The procedure results in a Lie-Poisson system, which can be considered a generalization of a Hamiltonian system. The Lie-Poisson reduction reveals the manifold on which the dynamics naturally occur. We can interpret the singularities in previously studied reductions as arising from poorly chosen projections of the solution manifold.

Our earlier work~\cite{Anurag:2024} contains a history of this well-studied problem, which we summarize here. Gröbli initiated the study of three interacting vortices in his 1877 doctoral thesis~\cite{Grobli:1877}.%
\footnote{We have translated Gröbli's thesis into English, and posted the results on arXiv.org~\cite{Goodman:2024}.}
Helmholtz had earlier derived the equations of motion for a system of \(N\) point vortices and solved them for \(N\le 2\), and Kirchhoff put them in Hamiltonian form~\cite{Helmholtz:1858, Kirchhoff:1876}. The solutions with \(N=2\) are simple---the vortices' speeds and the distances between them remain constant for all time. Gröbli was the first to systematically study the problem with \(N=3\), the smallest system for which these quantities evolve in time~\cite{Aref:1992}.

Gröbli's work predates the development of the qualitative theory of differential equations---the Encyclopedia of Mathematics credits Poincaré and Lyapunov with initiating this subject in the 1880s and 1890s~\cite{qualde-encyc-maths}. The nineteenth-century paradigm was to reduce the equations to quadratures, i.e., to separate variables and integrate to produce exact solution formulas in terms of elementary functions and the classical special functions. Gröbli found solutions by quadrature for the three-vortex problem and for specific symmetric initial configurations of four or more vortices.

Given three vortices, with circulations \(2\pi \Gamma_j\) and positions \((x_j,y_j)\), for \(j=1,2,3\), Gröbli defined a closed system for the squared pairwise vortex distances
\begin{equation}\label{l_eqn}
\dv{t} \begin{pmatrix} \ell_{23}^2 \\ \ell_{31}^2 \\ \ell_{12}^2 \end{pmatrix} =
4 \sigma   A 
\begin{pmatrix} 
   \Gamma_1 \left(\ell_{12}^{-2} -\ell_{31}^{-2} \right)\\
   \Gamma_2 \left(\ell_{23}^{-2} -\ell_{12}^{-2} \right)\\
   \Gamma_3 \left(\ell_{31}^{-2} -\ell_{23}^{-2} \right)
\end{pmatrix}.
\end{equation}
Here \(\ell_{ij}\) is the distance between vortices \(i\) and \(j\), \(A\) is the area of the triangle formed by the vortices, and \(\sigma=\pm 1\) depending on whether the three vortices are arranged clockwise or counterclockwise. The area \(A\) can be determined from the distances \(\ell_{i,j}\) using Heron's formula, thus closing the system. 

Two modern ideas are notably absent from Gröbli's work: phase-plane reasoning and Hamiltonian reduction. Poincaré was the first to make systematic use of phase portraits; Book II of his New Methods of Celestial Mechanics contains several examples~\cite{Poincare1893}. While Gröbli draws several accurate illustrations of vortex trajectories, none of these is a phase portrait, and reasoning with phase portraits is absent.  Second, Gröbli mentions Kirchhoff's Hamiltonian interpretation of Helmholtz's equations, but does not perform canonical reductions. That is, he does not use canonical transformations to preserve the system's Hamiltonian structure while reducing its dimension.

Several authors independently found ways to apply phase-plane reasoning to Gröbli's reduced system~\cite{Aref:1979, Synge.1949, Novikov.1975}, noting that the conserved angular impulse can be used to eliminate one variable from system~\eqref{l_eqn}. However, the \(\ell_{i,j}\) coordinates only describe a physical situation when they satisfy the triangle inequality. Thus, the accessible phase space is a proper subset \(\mathcal{P} \subset \Reals^2\). The triangle inequality holds as an equality on \(\partial\mathcal{P}\), where the three particles are collinear, and the area vanishes. Heron's formula has a square-root singularity on \(\partial\mathcal{P}\), so that the right-hand side of system~\eqref{l_eqn} is not differentiable there. This prevents the use of linearization, the most basic method for determining stability. Since working with these equations is difficult, researchers have derived many alternate reductions, each with its own strengths and weaknesses~\cite{Krishnamurthy:2018, Stremler:2021, Luo:2022, Borisov:1998, Bolsinov:1999}. 

Questions of continuation and bifurcation, i.e., how changes in parameters affect the phase space, are important when studying any dynamical system. In the vortex problem, the circulations of the vortices serve as the bifurcation parameters, and relative equilibria---solutions that appear as rigid motions when viewed in an appropriate moving reference frame---are the principal invariant objects of interest. Previous groups have not found system~\eqref{l_eqn} useful for this task and have tackled this question by deriving alternate forms of the equations~\cite{Conte:1979, Conte:2015, Tavantzis.1988, Aref:2009}. Two of these groups have drawn a bifurcation diagram that encapsulates how the number of relative equilibria and their stability change with the parameters; see Fig.~\ref{fig:trilinear}. The methods they describe are useful only for finding the relative equilibria and their stability, not for drawing phase portraits. 

Following that work, the body of this paper considers the case of nonvanishing total circulation. By rescaling time, and possibly changing its sign, we may, without loss of generality, consider the case
\begin{equation} \label{circulation_1}
    \Gamma_1 + \Gamma_2 + \Gamma_3 =1.
\end{equation}
The case of vanishing total circulation is considered briefly in~\ref{sec:vanishing}.

The preferred approach to studying Hamiltonian systems, which we follow here, exploits conservation laws and symmetries to obtain a reduced system that possesses a Hamiltonian structure. This often proceeds in stages~\cite{Marsden:2007}.

Our first stage uses Jacobi coordinates, which exploit the system's translation invariance; this step is unchanged from our prior work~\cite{Anurag:2024}. The second stage uses a momentum map to mod out the rotational invariance, following recent papers of Ohsawa~\cite{Ohsawa:2019, Ohsawa:2024, Ohsawa:2025}. Finally, we introduce coordinates on the symplectic leaves of the reduced system, using a conservation law to reduce the number of variables from four to three. Because the solutions are confined to a two-dimensional manifold in \(\mathbb{R}^3\), the dynamics are equivalent to a one-degree-of-freedom system.

The reduced equations are the simplest yet produced, so that the analysis will be the most straightforward. Compared with prior results, the reduced equations involve neither square roots nor trigonometric functions, and the phase portrait lacks the singularities or inaccessible regions seen elsewhere.

The most important consequence of this reduction is that it gives both types of information in a single representation. First, it lets us construct a phase portrait representation that is simple to work with algebraically and interpret. At the same time, it allows us to compute and explain the bifurcation diagram in Fig.~\ref{fig:trilinear}.

We can interpret our final coordinate system as describing how the triangle formed by the three vortices evolves in time, similar to Montgomery's shape-space description of the gravitational \(3\)-body problem~\cite{Montgomery:2015}. For the \(3\)-vortex problem, the shape space may be compact, as in the gravitational problem, or unbounded. Ohsawa has derived a similar shape-space representation of three-vortex motion~\cite[Fig. 4]{Ohsawa:2019}. The Jacobi coordinates and our final parameterization of the symplectic leaves help simplify this picture.

This paper is organized as follows:
In Sec.~\ref{sec:vortex_hamiltonian}, we introduce the point vortex model and its Hamiltonian formulation. We review the two reduction methods used in the next section.
Sec.~\ref{sec:reduction} applies these methods to arrive at the final form of the equations in the case of nonvanishing total circulation.
Sec.~\ref{sec:bifurcations} analyzes the reduced equations to describe the fixed points and their stability. This section uses discriminants and resultants computed in Mathematica to identify curves in parameter space along which bifurcations occur, summarized in Fig.~\ref{fig:trilinear}. 
In Sec.~\ref{sec:phase_portraits}, we plot the phase portrait for the dynamics at different points in the parameter space described by the figure, concentrating on the case \(\Gamma_1=\Gamma_2\). This extra symmetry simplifies the equations sufficiently to find all the equilibria in closed form and introduces additional symmetry into the phase portraits.
We conclude in Sec.~\ref{sec:conclusion} and discuss other problems to which our analysis should apply.

The paper includes three appendices. 
In the first,~\ref{sec:vanishing}, we consider the case of vanishing total circulation, where the Jacobi coordinates of Sec.~\ref{sec:reduction} break down. These results are equivalent to those derived by Rott and Aref, and are included for completeness.
In~\ref{sec:resultant}, we discuss the resultant and the discriminant, two functions that allow us to analyze the polynomials whose roots determine the equilibria and their stability.
\ref{sec:trilinear} describes the barycentric coordinates used to construct Fig.~\ref{fig:trilinear}.

\section{The point-vortex model and Hamiltonian background}
\label{sec:vortex_hamiltonian}
The motion of \(N\) vortices with positions \(\r_j = (x_j, y_j)\) and nonzero circulations \(2\pi\Gamma_j\) evolves according to
\begin{equation} \label{N_vortex_equations}
\dv{x_j}{t} = - \sum_{\substack{i=1 \\ i\neq j}}^N \frac{\Gamma_i(y_j-y_i)}{\norm{\r_j-\r_i}^2}, \quad
\dv{y_j}{t} =  \sum_{\substack{i=1 \\ i\neq j}}^N \frac{\Gamma_i(x_j-x_i)}{\norm{\r_j-\r_i}^2}.
\end{equation}
This can be written as a Hamiltonian system of equations
\begin{equation} \label{N_vortex_ham_eqns}
  \dv{x_j}{t}  = \frac{1}{\Gamma_j}\pdv{H}{y_j}, \quad 
  \dv{y_j}{t}  = \frac{-1}{\Gamma_j}\pdv{H}{x_j},
\end{equation}
where
\begin{equation} \label{N_vortex_hamiltonian}
    H(\r)=\frac{-1}{2} \sum_{i< j}^N  \Gamma_j \Gamma_j   \log \norm{ \r_i-\r_j }^2.
\end{equation} 
Helmholtz first derived system~\eqref{N_vortex_equations}, and Kirchhoff put it in Hamiltonian form~\cite{Helmholtz:1858, Kirchhoff:1876}. It is a standard topic in elementary courses in mathematical fluid mechanics~\cite{Chorin:1993} and is covered thoroughly in Newton's textbook~\cite{Newton.2001}.

We provide a brief overview of the Hamiltonian formulation of this system, following the presentation in~\cite{Ohsawa:2024, Ohsawa:2025}, which contains a more complete description. It is convenient to arrange the vortex positions into a vector \(\r =\parentheses{x_1,\ldots,x_N,y_1,\ldots,y_N}\in \Reals^{2N}\) and to identify this with a vector \(\z \in \Complex^N\). 
We define the \(N\times N\) matrix
\begin{equation}\label{DGamma}
\DGamma \defeq \diag\parentheses{\Gamma_1,\ldots,\Gamma_N}
\end{equation}
and the \((2N)\times(2N)\) Poisson matrix
\begin{equation}\label{JJ}
\JJ \defeq \begin{bmatrix}
0 & \DGamma \\ 
-\DGamma & 0
\end{bmatrix}.
\end{equation}
The function \(H(\z)\) defines a Hamiltonian vector field \(X_H\) on \(\Reals^{2N} \cong \Complex^N\) by
\begin{equation}\label{XH}
X_H(\z) \defeq \left(\JJ^\Tee\right)^{-1} D H(\z) = -\JJ^{-1} DH(\z),
\end{equation}
where 
\begin{equation}\label{JTinv}
\left(\JJ^\Tee\right)^{-1} = - \JJ^{-1}=
\begin{bmatrix}
0 & \DGamma^{-1} \\ 
-\DGamma^{-1} & 0
\end{bmatrix},
\end{equation}
and \(\DGamma^{-1} = \diag\parentheses{\Gamma_1^{-1},\ldots,\Gamma_N^{-1}}\).

This allows us to define a circulation-dependent Poisson bracket of two \(C^1\) functions from \(\Complex^N\) to \(\Complex\):
\begin{equation} \label{poisson_bracket}
\begin{split}
\poisson{F(\z)}{G(\z)} 
& = X_F(\z)^\Tee\JJ X_G(\z)\\
& = DF(\z)^\Tee \left(\JJ^\Tee\right)^{-1} DG(\z) \\
& = \sum_{j=1}^N \frac{1}{\Gamma_j} \left( \pdv{F}{x_j} \pdv{G}{y_j} - \pdv{F}{y_j} \pdv{G}{x_j} \right).
\end{split}
\end{equation}
One consequence of the Poisson bracket formulation is that, using system~\eqref{N_vortex_ham_eqns}, any function \(F(\z)\) evolves according to 
\begin{equation} \label{Fdot}
\dv{t} F(\z) = X_H(F) \equiv \poisson{F(\z)}{H(\z)},
\end{equation}
and, in particular
\begin{equation} \label{zdot}
\dv{\z}{t}  = X_H(\z) \equiv \poisson{\z}{H(\z)}.
\end{equation}

This system has three conservation laws
\begin{equation}\label{constants_of_motion}
    \vb{M}=M_x + \rmi M_y= \sum_{j=1}^N\Gamma_j z_j,  \qand
    \Theta =\sum_{j=1}^N \Gamma_j \abs{z_j}^2.
\end{equation}
The quantity \(\vb{M}\) is the \emph{linear impulse}, and \(\Theta\), the \emph{angular impulse}. In the case that \( \sum_{j=1}^N \Gamma_j \neq 0\), then 
\begin{equation}\label{center_of_vorticity}
\z_0 = \vb{M}/\sum_{j=1}^N \Gamma_j = \vb{M}/\gamma_1
\end{equation}
defines the conserved center of vorticity. 
For later reference, we define the three \emph{elementary symmetric polynomials} in three circulations,
\begin{equation}\label{gammadef}
\gamma_1 = \Gamma_1 + \Gamma_2 + \Gamma_3, \quad
\gamma_2 = \Gamma_1  \Gamma_2 + \Gamma_3  \Gamma_1 + \Gamma_2  \Gamma_3, \qand
\gamma_3 = \Gamma_1  \Gamma_2  \Gamma_3.
\end{equation}
Any symmetric polynomial can be written as a polynomial in the \(\gamma_j\).

The conservation laws arise due to the system's invariance under orientation-preserving rigid transformations. That is, the system has \(\SE(2) = \SO(2)\ltimes\Reals^2\) symmetry under the action
\begin{equation}\label{SE2symmetry}
\SE(2)\times\Complex^N \to \Complex^N ; \;
\parentheses{ \parentheses{e^{\rmi\psi},a}, \z} \mapsto e^{\rmi \psi} \z + a \vone,
\end{equation}
where \(
\vone = \parentheses{1,\ldots,1}^\Tee\). 

We define some terms that may be unfamiliar. \(\SE(2)\) is the special Euclidean group of orientation-preserving rigid deformations of the plane. Any such deformation is the composition of a rotation \(e^{\rmi \psi} \in \SO(2)\) followed by a translation by \(a \in \Reals^2\). The "\(\ltimes\)" symbol denotes that \(\SE(2)\) is the semidirect product of the rotation and translation groups in the plane.

The two changes of variables described below reduce the dimension of the system. The first, a Jacobi coordinate reduction, removes the translations by \(a\vone\). The second, which puts the equations in Lie-Poisson form, eliminates the rotations by \(e^{\rmi \psi}\). We prefer \emph{canonical} changes of variables, i.e., those that preserve the system's Hamiltonian structure. Doing this without introducing new singularities will require us to use a more general definition of Hamiltonian mechanics than has been common in the point-vortex literature.

\subsection{Jacobi Coordinate Reduction}
The Jacobi coordinate reduction is a canonical change of variables that allows us to systematically introduce the center of vorticity \(\z_0\), defined in~\eqref{center_of_vorticity}, as a variable. Because \(\z_0\) does not vary, this reduces the number of degrees of freedom by one. The procedure constructs a change of the coordinates \(z_j\), and also introduces new circulation variables \(\tilde{\Gamma}_j\) into the Poisson bracket~\eqref{poisson_bracket}. It proceeds iteratively, replacing two conjugate pairs of variables at a time, repeating the procedure \(N-1\) times until all the variables have been replaced. Consequently, we can define the change of variables for the case \(N=2\) and iterate. 

Consider a Hamiltonian system of the form~\eqref{N_vortex_hamiltonian} with \(N=2\) and \(\Gamma_1 + \Gamma_2 \neq 0\). Define the new "position" variables and \emph{reduced circulations}
\begin{equation} \label{Jacobi1}
\begin{aligned}
 \tilde{z}_1 &= z_1 - z_2; & \tilde{\Gamma}_1 & = \frac{\Gamma_1 \Gamma_2}{\Gamma_1 + \Gamma_2}; \\
 \tilde{z}_2 & = \frac{\Gamma_1 z_1 + \Gamma_2 z_2}{\Gamma_1+\Gamma_2}; & \tilde{\Gamma}_2 & = \Gamma_1 + \Gamma_2,
\end{aligned}
\end{equation}
Letting 
\[
\tilde{H}\parentheses{\tilde{z}_1,\tilde{z}_2} = 
H\parentheses{z_1(\tilde{z}_1,\tilde{z}_2),z_2(\tilde{z}_1,\tilde{z}_2)},
\]
the system evolves under a system of the form~\eqref{zdot}, but in the "tilde" variables. The variable \(\tilde{z}_1\) is the displacement vector from \(z_2\) to \(z_1\), while \(\tilde{z}_2\) coincides with the pair's center of vorticity. 

If we were only interested in the case \(N=2\), the procedure would end here. From this calculation, we recover the elementary result that the two vortices orbit their mutual center of vorticity in a circle at a constant angular velocity.

\begin{figure}[htbp] %
   \centering
   \includegraphics[width=0.42\textwidth]{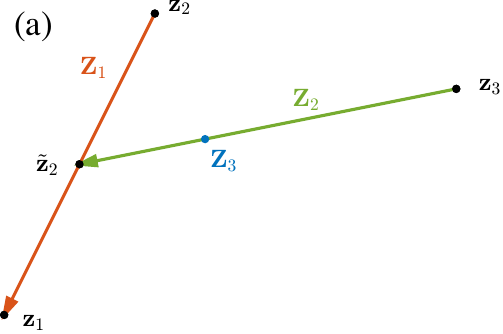}
   \hspace{0.1\textwidth}
   \includegraphics[width=0.30\textwidth]{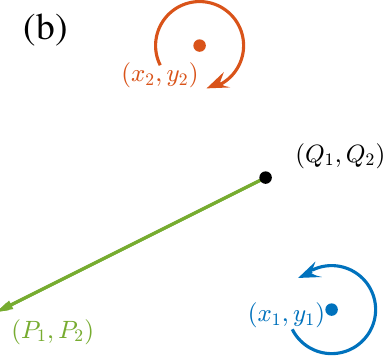}
   \caption{\textbf{(a)} Jacobi coordinates for three vortices with nonzero total circulation. The coordinate \(Z_1\) is the vector from \(z_2\) to \(z_1\), \(Z_2\) is the vector from \(z_3\) to the \(\tilde{\z}_2\) (the center of vorticity of the first two particles, and the transformed variable \(Z_3\) is the conserved center of vorticity. \textbf{(b)} The alternative to Jacobi coordinates for two point-vortices with \(-\Gamma_2=\Gamma_1>0\). }
   \label{fig:Jacobi}
\end{figure}

When \(N\ge3\), the transformation to Jacobi coordinates is defined iteratively: the change of variables~\eqref{Jacobi1} is the first step, while the variables with subscripts \(3\ldots N\) are left unchanged. Next, combine \(\tilde{z}_2\) with \(z_3\), using formulas analogous to~\eqref{Jacobi1}, and continue with coordinates labeled \(k\) and \((k+1)\) until all the position variables and circulations have been transformed in this way. Fig.~\ref{fig:Jacobi}(a) demonstrates the case \(N=3\). The change of variables will be undefined at step \(k\) if \(\sum_{j=1}^{k+1}\Gamma_j=0\). If the total circulation is nonzero, then a particle ordering always exists that guarantees this problem will never occur. 

The procedure yields a new position vector \(\Z=(Z_1,\ldots,Z_{N-1})\) and reduced circulation vector \((\kappa_1,\ldots,\kappa_{N-1})\). We have not included the \(Z_N\) and \(\kappa_N\) in these vectors since the former corresponds to the conserved center of vorticity. The Hamiltonian is independent of \(Z_N\), so the system has been reduced by one degree of freedom. 

The change to Jacobi coordinates depends on the order of the labels assigned to the point vortices. All such reductions lead to equivalent systems, but some may be more convenient to work with than others. When \(\gamma_1=0\), the final application of Eq.~\eqref{Jacobi1} is undefined because the two fractions appearing in that equation have vanishing denominators. We discuss this case in~\ref{sec:vanishing}.

\subsection{Reduction to Lie-Poisson form}
\label{sec:lie-poisson}

The reduction below is possible because the  Hamiltonian~\eqref{N_vortex_hamiltonian} and the Poisson bracket~\eqref{poisson_bracket} are both invariant under the \(\SS^1\) action on \(\Complex^N\), 
\begin{equation}\label{S1symmetry}
\SS^1\times\Complex^N \to \Complex^N ; \;
\parentheses{ e^{\rmi\psi}, \z} \mapsto e^{\rmi \psi} \z.
\end{equation}
We introduce a change of variables to a new variable \(\mu\) whose form may appear surprising at first glance. The simplification is not yet apparent, but it will emerge in the subsequent analysis. 
We define the function 
\begin{equation} \label{z_to_mu}
\J: \Complex^N \to \ufrak(N); \quad 
\z \mapsto \mu 
 \coloneqq
\rmi \z \z^* ,
\end{equation}
where the asterisk represents the Hermitian transpose. \(\J\) sends complex vectors to complex skew-Hermitian matrices, i.e., those satisfying \(\mu^*=-\mu\). Following Marsden and Ratiu~\cite{Marsden:1999}, Ohsawa notes that \(\J\) is a momentum map~\cite[Sec 2.3]{Ohsawa:2025}.
For \(N=2\), we write
\begin{equation} \label{mu2by2}
\mu 
= \rm i \begin{bmatrix}
\mu_1 & \mu_3 + \rmi \mu_4 \\
 \mu_3 - \rmi \mu_4 & \mu_2
\end{bmatrix} 
 =
\rmi \begin{bmatrix}
|z_1|^2 & z_1^{\phantom{*}} z_2^* \\
z_{2}^{\phantom{*}} z_1^* & |z_2|^2
\end{bmatrix}.  
\end{equation}

Since the Hamiltonian is invariant to \(\SS^1\) transformations of the form~\eqref{S1symmetry}, it can be defined entirely in terms of elements of the matrix \(\mu\). Thus we may write the Hamiltonian as \(h(\mu)\), where
\[
H = h\circ\J.
\]
The \(\mu\) coordinates can be shown to satisfy the evolution equation
\begin{equation} \label{mudot}
\dot{\mu} = \DGamma^{-1} \fdv{h}{\mu} \mu - \mu \fdv{h}{\mu} \DGamma^{-1},
\end{equation} 
where the functional derivative is
\begin{equation} \label{dfdmu2}
\fdv{f}{\mu} = 
\rmi
\begin{bmatrix}
2 \pdv{f}{\mu_1} &&  \pdv{f}{\mu_3} + \rmi \pdv{f}{\mu_4}\\
\pdv{f}{\mu_3} - \rmi \pdv{f}{\mu_4} && 2\pdv{f}{\mu_2},
\end{bmatrix}
\end{equation}
for \(N=2\) and takes a similar form for general \(N\) as can be demonstrated using Eq.~\eqref{dfdmu} below.

This can be interpreted as a Hamiltonian system because this change of variables is a Poisson transformation, which we explain below.  Equations~\eqref{z_to_mu}--\eqref{dfdmu2} are sufficient for a reader who wants to work through the calculations for themselves. We briefly describe their derivation and meaning, and refer to Ohsawa~\cite{Ohsawa:2019, Ohsawa:2024, Ohsawa:2025} for more details and an application of this reduction to several stability problems in vortex dynamics, and to  Marsden and Ratiu for the necessary background in geometric mechanics~\cite{Marsden:1999}. 

First, define \(\ufrak(N)_\Gamma\) as the vector space of skew-Hermitian \(N\times N\) matrices equipped with the bracket
\[
\\textbf{brackets}{\xi,\eta}_\Gamma = \xi \DGamma^{-1} \eta - \eta \DGamma^{-1} \xi, \qfor \xi, \eta \in \ufrak(N)_\Gamma,
\]
an antisymmetric bilinear operator that satisfies the Jacobi identity. Therefore, \(\ufrak(N)_\Gamma\) is a Lie algebra.

Next, define the inner product
\begin{equation} \label{innerproduct}
\pair{\xi}{\eta} = \frac{1}{2} \tr\parentheses{\xi^* \eta}.
\end{equation}
Now, \(\ufrak(N)_\Gamma^*\), the dual of \(\ufrak(N)_\Gamma\), is the set of linear operators on \(\ufrak(N)_\Gamma\), and can be identified bijectively with \(\ufrak(N)_\Gamma\): if \(\alpha \in \ufrak(N)_\Gamma^*\), there exists a unique \(\alpha^\sharp \in \ufrak(N)_\Gamma\) such that
\[
\alpha(\mu) = \pair{\alpha^\sharp}{\mu}
\]
for all \(\mu \in \ufrak(N)_\Gamma\). 

Inner product~\eqref{innerproduct} lets us define a Poisson bracket on functions \(f\) and \(g\) from \(\ufrak(N)_\Gamma^* \cong\ufrak(N)_\Gamma\) to \(\Reals\),
\begin{equation} \label{poissonGamma}
\poisson{f}{g}_\Gamma(\mu) \defeq \pair*{\mu}{ \brackets{\fdv{f}{\mu},\fdv{g}{\mu}}_\Gamma}.
\end{equation}
Using these definitions, we return to our definition of the momentum map in Eq.~\eqref{z_to_mu} and note that the codomain is properly defined as \(\ufrak(N)_\Gamma^*\), because, as discussed in~\cite{Marsden:1999}, this Poisson bracket, actually a \emph{Lie-Poisson bracket}, is defined on the dual of a Lie algebra.

The functional derivatives are defined by enforcing that for any \(\mu,\nu \in \ufrak(N)_\Gamma^*\)
\begin{equation} \label{dfdmu}
\pair*{\nu}{\fdv{f}{\mu}} 
= \frac{1}{2} \tr \parentheses{ \nu^* \fdv{f}{\mu}} 
= \left. \dv{s} \right|_{s=0} f(\mu + s \nu) 
= \sum \left(\fdv{f}{\mu}\right)\!{}_{\!j} \nu_j.
\end{equation}
Equating the trace-valued inner product and the explicit sum for \(N=2\) yields Eq.~\eqref{dfdmu2}.

Putting this all together, we can now say precisely how \(\J\) is canonical. It is a \emph{Poisson map }with respect to the Poisson brackets~\eqref{poisson_bracket} and~\eqref{poissonGamma}, meaning that
\[
\poisson{f\circ\J}{g\circ\J} = \poisson{f}{g}_\Gamma \circ \J
\]
for arbitrary smooth functions \(f\) and \(g\). Therefore, any smooth function \(f(\mu)\) evolves under
\[
\dv{t} f(\mu) = \poisson{f}{h}_\Gamma,
\]
which has the same functional form as Eq.~\eqref{zdot}.
Applying this to the components of \(\mu\) gives the evolution equation~\eqref{mudot}.

\begin{figure}[htbp] %
   \centering
   \includegraphics[width=2in]{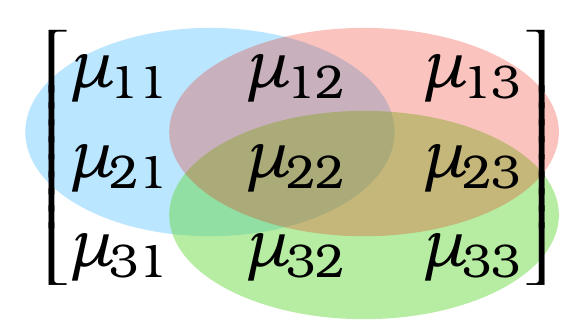} 
   \caption{The three independent subdeterminants that must vanish in the case \(N=3\).}
   \label{fig:determinant_condition}
\end{figure}

The final question is in what sense an equation for the \(N^2\) real-valued components of the skew adjoint matrix \(\mu\) is simpler than an equation for the \(2N\) complex-valued components of \(\z\) or for previously derived reduced models. 
The first is that the simplest states in the \(N\)-vortex model are relative fixed points: they rigidly translate or rotate while maintaining a constant shape. These become fixed points in the reduced system.
Second, while many changes of variables introduce singularities, Eq.~\eqref{z_to_mu} does not; for example, Gröbli's equations become singular whenever the three vortices are collinear. Finally, since \(\mu\) is defined as the exterior product of two vectors, it has rank one. Therefore, all the rows except the first are redundant. A practical consequence is that each \(2\times2\) subdeterminant of \(\mu\) must vanish. For matrix~\eqref{mu2by2}, this results in the single real equation
\begin{equation}\label{rankone2by2}
\mu_1 \mu_2 - \mu_3^2 - \mu_4^2 = 0.
\end{equation}
For systems of \(N\) equations, this results in \((N-1)^2\) real conditions as follows. The functions
\[
\delta_{ij} = \det \begin{pmatrix} \mu_{ij} & \mu_{i,j+1} \\ \mu_{i+1,j} & \mu_{i+1,j+1} \end{pmatrix} 
\qfor 1\le i \le N-1 \qand i\le j \le N-1,
\] 
must vanish. For \(j=i\), this provides a real-valued constraint because \(\mu^*=-\mu\), and for \(j>i\), the constraint is complex valued, so its real and imaginary parts both vanish. Subtracting the number of constraints from \(N^2\) leaves \(2(N-1)\) free variables, reducing the degrees of freedom by one. The required subdeterminants for a system of three equations are illustrated in Fig.~\ref{fig:determinant_condition}.

\section{The reduced equations for \(N=3\)}
\label{sec:reduction}

We now specialize to the case \(N=3\) in two stages, first, applying Eq.~\eqref{Jacobi1}, to vortices 1 and 2. We can always label the vortices such that \(\Gamma_1+\Gamma_2 \neq 0\), and the reduction is valid. Next, we apply Eq.~\eqref{Jacobi1} to the variables \(\tilde{z}_2\) and \(\tilde{z}_3\). Calling the transformed position coordinates \(Z_j\) and the virtual circulations \(\kappa_j\), we find
\begin{equation}\label{jacobi_coordinates}
\begin{aligned}
    Z_1 & =  z_1 - z_2; & 
    \kappa_1 & = \frac{\Gamma_1  \Gamma_2}{\Gamma_1 + \Gamma_2}; \\
    Z_2 & = \frac{\Gamma_1z_1 + \Gamma_2 z_2}{\Gamma_1 + \Gamma_2} -z_3; 
    & \kappa_2 & =  \frac{(\Gamma_1 + \Gamma_2)\Gamma_3}{\Gamma_1 + \Gamma_2 + \Gamma_3}; \\
    Z_3 & = \frac{\Gamma_1 z_1 + \Gamma_2 z_2 + \Gamma_3 z_3}{\Gamma_1+\Gamma_2+\Gamma_3}; &  \kappa_3 &= \Gamma_1 + \Gamma_2 +\Gamma_3.
\end{aligned}
\end{equation}
These are simplified by using Eq.~\eqref{circulation_1}.

The variable \(Z_3\) is the conserved center of vorticity. Without loss of generality, we may choose coordinates so that \(Z_3=0\), yielding the  Hamiltonian
\begin{equation} \label{generalhamiltonian}
  H = \frac{-\Gamma_1 \Gamma_2}{2} \log{\abs{Z_1}^2} -
       \frac{\Gamma_2  \Gamma_3}{2} \log{\abs{Z_2 - \frac{\kappa_1}{\Gamma_2}Z_1}^2} 
      - \frac{\Gamma_1 \Gamma_3}{2} \log{\abs{Z_2+ \frac{\kappa_1}{\Gamma_1}Z_1}^2},
\end{equation}
and the conserved angular impulse 
\begin{equation} \label{generalangular}
  \Theta = \kappa_1 \abs{Z_1}^2 + \kappa_2 \abs{Z_2}^2.
\end{equation}

The reduced Poisson bracket and evolution equations depend on the \emph{virtual positions} \(Z_1\) and \(Z_2\) and the \emph{virtual circulations} \(\kappa_1\) and \(\kappa_2\) in the same way that equations~\eqref{N_vortex_ham_eqns} and~\eqref{poisson_bracket} depend on the physical positions \(z_j\) and circulations \(\Gamma_j\). 

Applying the momentum-map reduction of Sec.~\ref{sec:lie-poisson} to the reduced system in \(Z_j\) and \(\kappa_j\) with \(N=2\) yields
\begin{equation}
  \label{mudef}
  \mu =
  \rmi
  \begin{bmatrix}
    \mu_{1} & \mu_{3} + \rmi \mu_{4} \\
    \mu_{3} - \rmi \mu_{4} & \mu_{2}
  \end{bmatrix}
  \defeq \rmi Z Z^{*} = \rmi
  \begin{bmatrix}
    |Z_{1}|^{2} & Z_{1} Z_{2}^{*} \smallskip\\
    Z_{2} Z_{1}^{*} &  |Z_{2}|^{2}
  \end{bmatrix},
\end{equation}
Writing \(Z_j = R_j e^{\rmi \phi_j}\) and \(\phi = \phi_1-\phi_2\) provides a geometric interpretation of the \(\mu_j\) variables:
\begin{equation*}
  \mu_{1} = R_{1}^{2},
  \qquad
  \mu_{2} = R_{2}^{2},
  \qquad
  \mu_{3} = R_{1} R_{2} \cos\phi,
  \qquad
  \mu_{4} = R_{1} R_{2} \sin\phi.
\end{equation*}
In particular, \(\mu_4\) is the signed area of the triangle formed by the vectors \(Z_1\) and \(Z_2\).
In these variables, the angular impulse~\eqref{generalangular} is 
\begin{equation} \label{ThetaMu}
\Theta = \kappa_1 \mu_1 + \kappa_2 \mu_2.
\end{equation}

System~\eqref{mudot} describes the evolution of a vector \(\mu \in\Reals^4\).
Since, by equations~\eqref{muXYZ} and~\eqref{hXYZ}, \(\pdv{h}{\mu_4}=0\), it evolves as
\begin{equation} \label{mudot_vec}
\dv{\mu}{t} =
\begin{bmatrix}
\frac{2}{\kappa_1} \pdv{h}{\mu_3} \mu_4 \medskip \\
 -\frac{2}{\kappa_2} \pdv{h}{\mu_3} \mu_4 \medskip \\
 2  \left(\frac{\pdv{h}{\mu_2}}{\kappa_2}-\frac{\pdv{h}{\mu_1}}{\kappa_1}\right)\mu_4 \medskip\\ 
 2 \mu_3 \left(\frac{\pdv{h}{\mu_1}}{\kappa_1}-\frac{\pdv{h}{\mu_2}}{\kappa_2}\right)
 + \pdv{h}{\mu_3} \left(\frac{\mu_2}{\kappa_1}- \frac{\mu_1}{\kappa_2}\right)
\end{bmatrix}.
\end{equation}

This system conserves three first integrals, the Hamiltonian \(h(\mu)\) as well as  
\[
 C_{1}(\mu)=\mu_{1} \mu_{2}-\mu_{3}^{2}-\mu_{4}^{2}, \qand
 C_{2}(\mu) = \Theta =\kappa_{1} \mu_{1}+\kappa_{2} \mu_{2},
\]
which are Casimirs of the Lie-Poisson bracket~\eqref{poissonGamma}. Their joint level sets are the symplectic leaves. The rank-one condition~\eqref{rankone2by2} implies that \(C_1=0\) on the physically relevant leaf. 

We therefore introduce coordinates parameterizing the symplectic leaves.
\begin{equation}\label{muXYZ}
X = \mu_3, \ Y = \mu_4, \qand
Z \defeq \kappa_1 \mu_1 - \kappa_2 \mu_2, 
\end{equation}
which render the rank-one condition~\eqref{rankone2by2} as
\begin{equation}\label{rankoneXYZ}
\Theta^2 = Z^2 + 4 \kappa_1 \kappa_2 \left(X^2+Y^2\right) 
= Z^2 + \frac{4\Gamma_1\Gamma_2\Gamma_3}{\Gamma_1+\Gamma_2+\Gamma_3}\left(X^2+Y^2\right).
\end{equation}
As \(\dv{\Theta}{t}=0\), trajectories are constrained to a quadric surface in \((X, Y, Z)\) space. Considering Eq.~\eqref{circulation_1}, this surface is an ellipsoid if one or three circulations are positive and one sheet of a two-sheeted hyperboloid if two are positive. In the ellipsoidal case, \(\Theta\), and \(Z\) may take either sign. When \(\kappa_1 \kappa_2<0\), \(\Theta\) may take either sign, and which sheet of the hyperboloid is physically relevant depends on both \(\Theta\) and the circulations.

Using the notation~\eqref{gammadef}, we may write the Hamiltonian as
\begin{equation}\label{hXYZ}
\begin{split}
h(X,Y,Z;\Theta) = &-\frac{\Gamma_1 \Gamma_2}{2}  \log{(Z+\Theta)^2}  \\
&-\frac{\Gamma_1 \Gamma_3}{2}  
\log{\left(\phantom{-}4\gamma_3 X + (\gamma_1 \Gamma_1 - \Gamma_2\Gamma_3)Z + \left(\Gamma_1+\Gamma_2\right) \left(\Gamma_1+\Gamma_3\right) \Theta \right)^2}\\
&-\frac{\Gamma_2 \Gamma_3}{2}  
\log{\left(-4\gamma_3 X + (\gamma_1 \Gamma_2 - \Gamma_1\Gamma_3)Z + \left(\Gamma_1+\Gamma_2\right) \left(\Gamma_2+\Gamma_3\right) \Theta \right)^2}
,
\end{split}
\end{equation}
yielding evolution equations, again, using \(\pdv{h}{Y}=0\),
\begin{equation}\label{XYZdot}
\dv{t} \mqty[X\\Y\\Z] 
\defeq \mathbf{F}(X,Y,Z)
=
\mqty[-4 h_{Z} Y \medskip \\
 4 h_{Z} X -\frac{1}{\kappa_1\kappa_2} h_{X}Z \medskip \\
4 h_{X} Y].
\end{equation}
This formula is valid regardless of \(\sign{\kappa_2}\) and thus unifies the two systems of equations derived in~\cite{Anurag:2024}. 

The Hamiltonian~\eqref{hXYZ} is singular when two vortices occupy the same location in space, where the argument of the corresponding logarithmic term vanishes. The three singularities must solve both \(Y=0\) and Eq.~\eqref{rankoneXYZ}. We list them for later reference:
\begin{equation} \label{XZ_singularities}
\begin{split}
\Singularity_{12} &= \parentheses{0,0,-\Theta}, \\
\Singularity_{13} &= \parentheses{\frac{-\gamma_1\Theta}{(\Gamma_1+\Gamma_2)(\Gamma_1+\Gamma_3)},0,
\frac{(\gamma_1 \Gamma_1 - \Gamma_2\Gamma_3)\Theta} {(\Gamma_1+\Gamma_2)(\Gamma_1+\Gamma_3)}},\\
\Singularity_{23} &= \parentheses{\frac{\gamma_1\Theta}{(\Gamma_1+\Gamma_2)(\Gamma_2+\Gamma_3)},0,
\frac{(\gamma_1 \Gamma_2 - \Gamma_1\Gamma_3)\Theta} {(\Gamma_1+\Gamma_2)(\Gamma_2+\Gamma_3)}}.
\end{split}
\end{equation}
The subscripts indicate which two vortices' positions coincide at the singularities. \(\Singularity_{13}\) and \(\Singularity_{23}\) diverge as  \(\Gamma_i+\Gamma_j \to 0\) for any pair of circulations.

The  singularities~\eqref{XZ_singularities} scale linearly with \(\Theta\), as do the equilibria~\eqref{equilateral_configs} and~\eqref{symmetric_equilibria}. As \(\Theta\to0\), they all approach the triple-collision singularity
\begin{equation} \label{tripSing}
\SingTriple = \parentheses{0,0,0}.
\end{equation}

\section{Local and global bifurcations}
\label{sec:bifurcations}

\begin{figure}[htb] %
   \centering
   \includegraphics[width=4in]{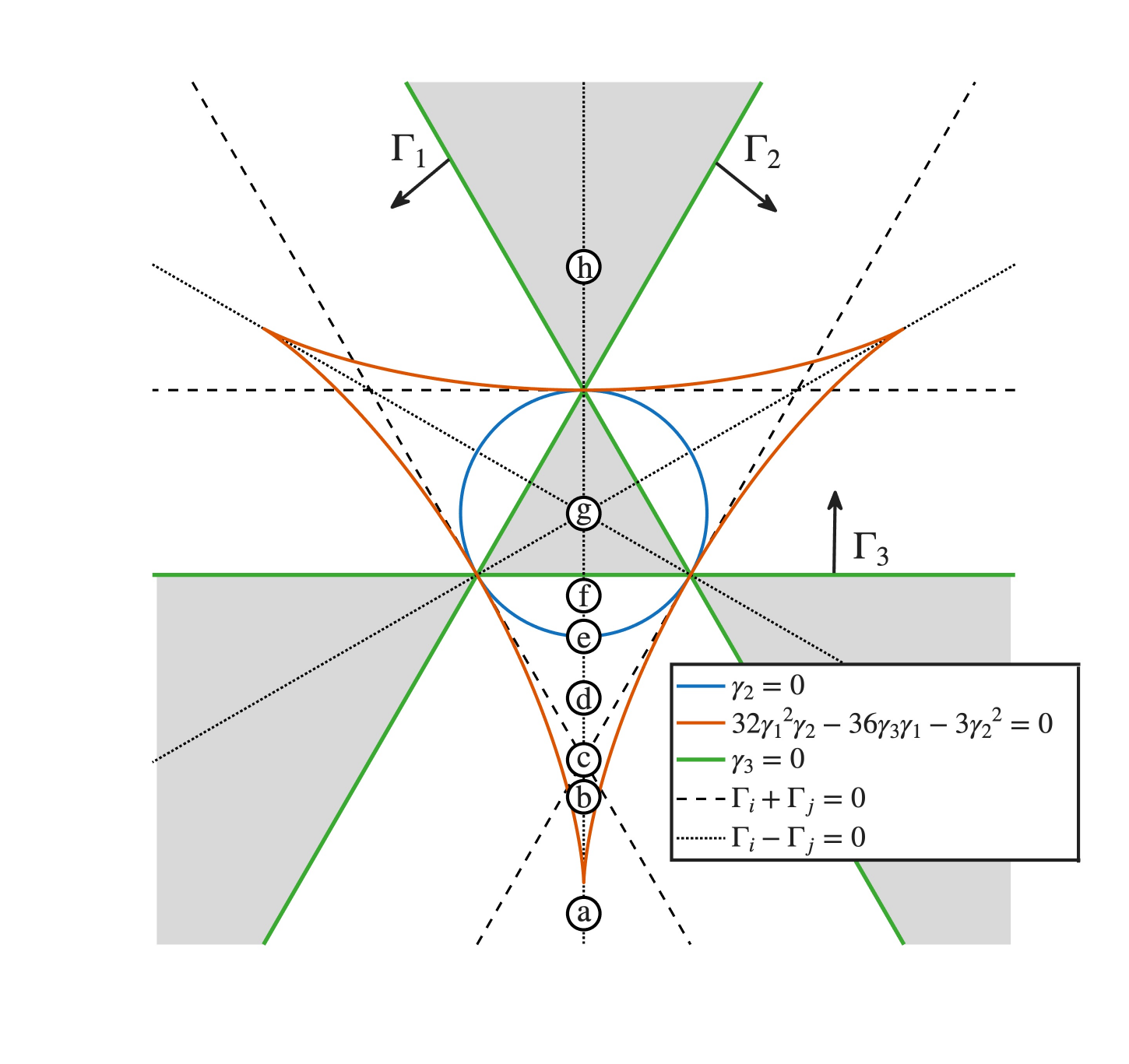} 
   \caption{The diagram identifying regimes of phase space behavior, defined by using \((\Gamma_1,\Gamma_2,\Gamma_3)\) as barycentric coordinates for the plane following Eq.~\eqref{barycenter_explicit}. The three arrows point to the direction in which the indicated circulation is positive, so that each vanishes on the line from which the arrow emerges. In the shaded region, the quadric surface defining the symplectic leaf in Eq.~\eqref{rankoneXYZ} is an ellipsoid, and in the unshaded region, it is a hyperboloid. The labels \textbf{a}-\textbf{h} indicate parameter values for which we draw phase planes in the next section.}
   \label{fig:trilinear}
\end{figure}

Central to this work is a bifurcation diagram in Fig.~\ref{fig:trilinear}, introduced by Conte in his 1979 Thèse d'État and revisited by him years later, as well as by Aref~\cite{Conte:1979, Conte:2015, Aref:2009}. It represents the space of circulations as \emph{barycentric coordinates}. Barycentric coordinates provide a bijective map between \(\Reals^2\) and the plane defined by~\eqref{circulation_1} as follows. Choose three noncollinear points \((\sfx_j,\sfy_j)\) in the plane. Then any point in \(\Reals^2\) can be written
\[
(\sfx,\sfy) = \sum_{j=1}^3 \Gamma_j\cdot(\sfx_j,\sfy_j).
\]
We explain this in more detail in~\ref{sec:trilinear}.
The reduced Hamiltonian evolution equations specified by Eqs.~\eqref{hXYZ} and~\eqref{XYZdot} allow us to analyze this diagram in an especially simple manner. 

The symplectic leaf is ellipsoidal if \(\Gamma_1 \Gamma_2 \Gamma_3>0\), corresponding to the shaded regions in the figure. It is a hyperboloid in the unshaded areas. Relative fixed points of system~\eqref{N_vortex_equations} are fixed points of system~\eqref{XYZdot}. It is well known that the three-vortex system has two families of relative fixed points: equilateral triangular and collinear configurations. We consider these in turn after discussing their stability.

Determining the linear stability of an equilibrium requires linearizing the system~\eqref{XYZdot} around that point. 
\begin{equation}\label{jacobian}
D\mathbf{F} = 
\begin{bmatrix}
 -4 h_{XZ} Y & -4 h_{Z} & -4 h_{ZZ} Y \\
 \frac{-1}{\kappa_1 \kappa_2}h_{XX} Z+4 \left(h_{XZ} X+ h_{Z} \right) & 0 & \frac{-1}{\kappa_1 \kappa_2}\parentheses{h_{X}+h_{XZ} Z}+4 h_{ZZ} X \\
 4 h_{XX} Y & 4 h_{X} & 4 h_{XZ} Y \\
\end{bmatrix},
\end{equation}
where subscripts denote partial derivatives.
Its determinant is
\begin{multline} \label{detDF}
   \det D\mathbf{F} = 
-16 Y \Bigg( \frac{h_X}{\kappa_1 \kappa_2} \left(h_X h_{XZ}-h_{XX} (h_{Z}+h_{ZZ} Z)+h_{XZ}^2 Z\right)\\+4 h_{Z} \left(h_{ZZ} (h_X+h_{XX} X) -h_{XZ}^2 X -h_{XZ} h_{Z}\right)\Bigg),
\end{multline}
and its trace is zero.

\subsection{Equilateral configurations}\label{sec:equilateral_equilibria}
When \(Y=0\), the three vortices are collinear, so \(Y\neq0\)  for the equilateral configurations. These satisfy \(h_{X}=h_{Z}=0\) and the rank-one condition~\eqref{rankoneXYZ}. We find two solutions
\begin{equation} \label{equilateral_configs}
\EqTriPM
= 
\frac{\Theta}{\gamma_2} 
\qty(
\frac{\left(\Gamma_1-\Gamma_2\right) \gamma_1}{2 \left(\Gamma_1+\Gamma_2\right)}, \;
\pm \frac{\sqrt{3} \gamma_1}{2}  , \;
\Gamma_1 \Gamma_2-\frac{\left(\Gamma_1^2+\Gamma_2^2\right) \Gamma_3}{\Gamma_1+\Gamma_2}).
\end{equation}
These fixed points diverge as \(\gamma_2 \to 0\) or \(\Gamma_1+\Gamma_2\to 0\). The equation \(\gamma_2=0\) parameterizes the circle circumscribing the central shaded triangle in Fig.~\ref{fig:trilinear}.

Linearizing about \(\EqTriPM\),  Jacobian~\eqref{jacobian} has characteristic polynomial
\begin{equation} \label{charPolyTri}
\chi_{\rm{tri}} (\lambda)  = - \lambda^3 - \frac{3 \gamma_2^3}{\Theta^2 \gamma_1^2}\lambda.
\end{equation}
The quadratic term vanishes because \(\tr D\mathbf{F}\equiv0\) and the constant term vanishes because the determinant~\eqref{detDF} vanishes when \(h_X=h_Z=0\). This always has one eigenvalue \(\lambda=0\). The other two eigenvalues indicate linear stability where \(\gamma_2>0\),inside the circle, and instability where \(\gamma_2<0\), outside the circle. 
 
\subsection{Collinear configurations}

On all collinear configurations, \(Y=0\), so that \(\dv{X}{t}\) and \(\dv{Z}{t}\) identically vanish by Eq.~\eqref{XYZdot}. Such equilibria satisfy both \(\dv{Y}{t}=0\) and the rank-one condition~\eqref{rankoneXYZ} with \(Y=0\). These are both polynomial equations
in \(X\) and \(Z\).%
\footnote{\(\dv{Y}{t}\) is a rational function, but we require that its numerator, a polynomial, vanishes.}
We may eliminate \(X\) from this system by computing the \emph{resultant}  of these two polynomials using Mathematica, as described in~\ref{sec:resultant}, yielding the condition
\begin{equation} \label{resultant}
\rho(Z;\Gamma_1,\Gamma_2,\Gamma_3,\Theta) = -64 \gamma_1 \gamma_3^2 \parentheses{Z-Z_{12}}\parentheses{Z-Z_{13}}\parentheses{Z-Z_{23}} p_3(Z;\Theta,\Gamma_1, \Gamma_2, \Gamma_3),
\end{equation}
where the \(Z_{ij}\) terms are the \(Z\)-components of singularities found in Eq.~\eqref{XZ_singularities} and \(p_3\) is cubic in \(\frac{Z}{\Theta}\) with the circulations appearing as parameters. The singularities cannot be solutions, so the system has one or three real roots, except at bifurcation points, where it has two. 

Two conditions indicate the parameter values where the number of solutions to \(p_3\) changes. First, when its leading coefficient $a_3$ vanishes, $p_3$ reduces to a quadratic, and so has two or fewer roots. As $a_3$ approaches zero, at least one root must diverge. The leading coefficient is
\begin{equation} \label{leading_coefficient}
\begin{split}
a_3&= \left(\Gamma_1+\Gamma_2\right)^2 \left(\Gamma_1+\Gamma_3\right) \left(\Gamma_2+\Gamma_3\right) \left(\Gamma_1 \Gamma_2+\Gamma_3 \Gamma_2+\Gamma_1 \Gamma_3\right)^2 \\
&=(\Gamma_1+\Gamma_2) \gamma_2^2 \left(\gamma_1 \gamma_2-\gamma_3\right).
\end{split}
\end{equation}
Thus, along the three lines \(\Gamma_i+\Gamma_j=0\), indicated as dashed lines in Fig.~\ref{fig:trilinear}, roots diverge to infinity. They also diverge where \(\gamma_2=0\), which corresponds to the unit circle in barycentric coordinates.

The number of solutions also changes at parameter values where the \emph{discriminant} of \(p_3(Z)\), which is proportional to the resultant of \(p_3\) and \(p_3'\), vanishes; this is also reviewed in~\ref{sec:resultant}. This quantity, computed in Mathematica, factors into terms of low order that we can interpret:
\begin{equation}\label{discriminant}
\Disc(p_3) = 64 \Theta^{6} 
   \parentheses{\Gamma_1-\Gamma_2}^2 \parentheses{\Gamma_1+\Gamma_2}^2 
   {\gamma_1}^2
   {\gamma_2}^2 
   {\gamma_3}^6 
   \left(32 \gamma_2 \gamma_1^2-36 \gamma_3 \gamma_1-3 \gamma_2^2\right).
\end{equation}
We review these terms below. Although we have assumed \(\gamma=1\), we leave it in this expression and those that follow.

\subsubsection*{The factor \(\Theta^{6}\)} When \(\Theta\neq0\), the evolution equations~\eqref{XYZdot} depend on \(X\), \(Y\), and \(Z\) only through terms of the form \(\frac{X}{\Theta}\), \(\frac{Y}{\Theta}\), and \(\frac{Z}{\Theta}\), so, for fixed values of the circulations, the phase space and fixed points scale with \(\Theta\) as long as \(\sign{(\Theta)}\) remains unchanged. This reflects the scale invariance of the \(N\)-vortex problem. 

\subsubsection*{The factors \(\parentheses{\Gamma_1-\Gamma_2}^2 \parentheses{\Gamma_1+\Gamma_2}^2\)}
Since the system~\eqref{N_vortex_equations} is invariant under permutations of the circulations, the bifurcation values must depend on the circulations only through the symmetric polynomials in Eq.~\eqref{gammadef}, so any factors that depend non-symmetrically on the circulations must be artifacts of the reduction method and must not indicate multiple roots. However, \((\Gamma_1+\Gamma_2)\) is a factor of the leading coefficient~\eqref{leading_coefficient}, as discussed above. Recall that the vanishing discriminant is necessary but insufficient for the two polynomials to vanish jointly.

\subsubsection*{The factor \(\left(32  {\gamma_1}^2\gamma_2-36 \gamma_3 \gamma_1-3 {\gamma_2}^2\right)\)} This factor is quartic in the circulations. In the figure, it vanishes along a curve with three arching segments and cusps at its three vertices, a well-known curve called the deltoid or Steiner's hypocycloid, defined as the image of a point on the circumference of a circle of radius one that rolls without slipping along the inside of a circle of radius three. We confirmed this by using Eq.~\eqref{barycenter_explicit}. Inside the deltoid, the system has three collinear equilibria; outside, it has only one. If the triple of circulations crosses the deltoid at a point other than the cusps, this arises due to a saddle-node bifurcation. The bifurcation is a pitchfork if the triple crosses through the cusp along the deltoid's symmetry axis. This is the generic behavior of the \emph{cusp catastrophe}~\cite[\S3.6]{Strogatz.2015} and is visible in the transition from Fig.~\ref{fig:point_a}(c) to Fig.~\ref{fig:point_b}(c).

\subsubsection*{The factor \(\gamma_2^2\)} This vanishes on the circle in Fig.~\ref{fig:trilinear} where the equilateral configurations change stability and the coefficient \(a_3\) vanishes, causing two collinear solutions to diverge. Since this term is square, the discriminant vanishes on the circle but does not change signs, and the number of solutions is the same inside and outside the circle, except at the three points where it is tangent to the deltoid.

\subsubsection*{The factor \({\gamma_3}^{6}\)} We may expand this factor as \({\Gamma_1}^{6} {\Gamma_2}^{6} {\Gamma_3}^{6}\). It vanishes along the three lines \(\Gamma_j=0\), but it does not change signs. These lines form the locus along which the geometry of the symplectic leaves bifurcates from ellipsoidal to hyperboloidal.

\subsubsection*{Linear Stability}
\label{sec:collinear_stability}
Plugging \(Y=0\) into the Jacobian~\eqref{jacobian}, the four entries at the corners of the matrix vanish identically, leaving a characteristic polynomial of the form
\begin{equation} \label{charPoly0}
\chi_{\mathrm{coll}}(\lambda) = -\lambda^3 - r(\Equilibrium) \lambda,
\end{equation}
where 
\begin{equation*} 
r(\Equilibrium) =
\frac{-4}{\kappa_1 \kappa_2} \left(h_{X}^2+h_{X} h_{XZ} Z-h_{XX} h_{Z} Z\right)-16 \left(-h_{X} h_{ZZ} X+h_{XZ} h_{Z} X+h_{Z}^2\right)
\end{equation*}
is a rational function in \(X\) and \(Z\) whose coefficients depend on \(\Theta\) and the circulations \(\Gamma_j\). An equilibrium \(\Equilibrium\) is stable when \(r(\Equilibrium)>0\).

Finding the collinear fixed points is equivalent to finding roots of the cubic polynomial \(p_3(Z)\) defined in Eq.~\eqref{resultant}. While this is formally solvable via Cardano's formula, the resulting expressions are too complicated to be practically useful. Therefore, we cannot simply evaluate \(r(\Equilibrium)\) at a collinear equilibrium to obtain a simple stability criterion as in Eq.~\eqref{charPolyTri} for the triangular equilibria.

Instead, we can search for values of the circulation at which the stability changes. At such equilibria, \(Y=0\), and the remaining components \(X\) and \(Z\) must satisfy three equations: condition~\eqref{rankoneXYZ}, \(\dv{Y}{t}=0\), and \(r(\Equilibrium)=0\). We may eliminate \(X\) and \(Z\) from this system by computing three resultants. First, choose two pairs of equations and eliminate \(X\) from each pair using resultants, producing two equations in \(Z\) alone. Then, compute a third resultant to eliminate \(Z\) from the first two, and  find the condition:
\begin{equation}\label{stability_resultant}
\Theta^{24} 
   p_{30}(\Gamma_1,\Gamma_2,\Gamma_3)
   \parentheses{\gamma_1 \gamma_2-\gamma_3}^2
   {\gamma_1}^{16}
   {\gamma_2}^{4} 
   {\gamma_3}^{24}
   \parentheses{32 \gamma_2 {\gamma_1}^2-36 \gamma_3 \gamma_1-3 {\gamma_2}^2} = 0.
\end{equation}
All the factors are symmetric in the three circulations except for \(p_{30}(\Gamma_1,\Gamma_2,\Gamma_3)\), a degree-30 homogeneous polynomial in the circulations. Like the factors \(\parentheses{\Gamma_1-\Gamma_2}^2 \parentheses{\Gamma_1+\Gamma_2}^2\) found in the discriminant formula above, these do not correspond to meaningful bifurcations, since they do not depend symmetrically on the circulations.

The three factors \({\gamma_1}^{16}
   {\gamma_2}^{4} 
   {\gamma_3}^{24}
   \parentheses{32 \gamma_2 {\gamma_1}^2-36 \gamma_3 \gamma_1-3 {\gamma_2}^2}\) vanish on the same sets as the discriminant~\eqref{discriminant}. The only new symmetric factor that arises in Eq.~\eqref{stability_resultant} and is not present in the discriminant is  
\[
 \parentheses{\gamma_1 \gamma_2-\gamma_3}^2=\parentheses{\Gamma_1+\Gamma_2}^2 \parentheses{\Gamma_1+\Gamma_3}^2 \parentheses{\Gamma_2+\Gamma_3}^2.
\] 
We recognize this factor from the coefficient~\eqref{leading_coefficient}. Because it is squared, it does not change sign and does not alter the stability.

\subsubsection*{The collinear equilibria when \(\Gamma_1=\Gamma_2\)}
\label{sec:Gamma23}
We can obtain more explicit results under the additional assumption \(\Gamma_1=\Gamma_2\). This occurs along the \(y\)-axis of Fig.~\ref{fig:trilinear}.  Assumption~\eqref{circulation_1} yields \(\Gamma_1=\Gamma_2 =\frac{1-\Gamma_3}{2}\). The cubic polynomial \(p_3(Z)\) that appears in the resultant~\eqref{resultant} factors as
\[
p_3(Z)= \frac{(1-\Gamma_3)^4}{64}\left(\left(\Gamma_3+1\right) \left(3 \Gamma_3+1\right) Z +\left(3 \Gamma_3^2+6 \Gamma_3-1\right) \Theta\right)^2(Z-\Theta),
\]
so that, after solving for \(X\), we have
\begin{equation}\label{symmetric_equilibria}
\begin{split}
\Equilibrium_1 &= \frac{\Theta}{(1+\Gamma_3)(1+3\Gamma_3)}
\parentheses{\phantom{-}\sqrt{\frac{5+3\Gamma_3}{1-\Gamma_3}},0,-3\Gamma_3^2-6\Gamma_3+1} ;\\
\Equilibrium_2 &= \frac{\Theta}{(1+\Gamma_3)(1+3\Gamma_3)}
\parentheses{-\sqrt{\frac{5+3\Gamma_3}{1-\Gamma_3}},0,-3\Gamma_3^2-6\Gamma_3+1}; \\
\Equilibrium_3 &= \parentheses{0,0,\Theta}.
\end{split}
\end{equation}

Clearly, the equilibrium \(\Equilibrium_3\) exists for all values of \(\Gamma_3\).
The equilibria \(\Equilibrium_1\) and \(\Equilibrium_2\) only exist when the term inside the square root is positive, that is, when \(\frac{-5}{3}<\Gamma_3<1\). 

On the line \(\Gamma_1=\Gamma_2\), there are three distinguished \(\Gamma_3\) intervals. 
\begin{itemize}
\item For \(\Gamma_3<0\), \(\kappa_1>0>\kappa_2\), so Eq.~\eqref{muXYZ} implies that \(Z\ge0\), making the symplectic leaf the upper sheet of the hyperboloid or upper half-cone. Depending on the sign of \(\Theta\), only the singularities or equilibria defined by Eqs.~\eqref{XZ_singularities},~\eqref{equilateral_configs}, and~\eqref{symmetric_equilibria} with \(Z>0\) lie on the sheet.
\item For \(0<\Gamma_3<1\), \(\kappa_1>0\) and \(\kappa_2>0\), so \(\Theta>0\) and \(Z\) may take either sign. The symplectic leaf is an ellipsoid.
\item For \(1<\Gamma_3\), \(\kappa_1<0\) and \(\kappa_2<0\), so \(\Theta<0\) and \(Z\) may take either sign.  The symplectic leaf is an ellipsoid.
\end{itemize}

We know from the argument following Eq.~\eqref{rankoneXYZ} that \(Z\ge0\) in the hyperboloid case. Therefore \(\Equilibrium_3\) lies on the leaves  \(\Theta\ge0\). The value \(\Gamma_3=\frac{-5}{3}\) occurs at the deltoid's cusp near the point \textbf{a} on the \(y\)-axis in Fig.~\ref{fig:trilinear} and \(\Gamma_3=1\) occurs at the upper edge of the deltoid where it intersects the \(y\)-axis. At this point, \(\Gamma_1=\Gamma_2=0,\) so the system degenerates to a one-vortex problem.

The function \(r(\Equilibrium)\) that determines stability in Eq.~\eqref{charPoly0} takes the values
\begin{equation*}
r(\Equilibrium_{1,2}) = \frac{\left(\Gamma_3-1\right)^2 \left(3 \Gamma_3+1\right)^3 \left(3 \Gamma_3+5\right) }{64 \Theta^2} \qand
r(\Equilibrium_3) = \frac{3 \left(\Gamma_3-1\right)^3 \left(3 \Gamma_3+5\right)  }{16 \Theta^2}.
\end{equation*}

The function \(r(\Equilibrium_1)\) is negative for \(\frac{-5}{3}<\Gamma_3<\frac{-1}{3}\), where \(\Equilibrium_1\) and \(\Equilibrium_2\) are centers and positive for \(\frac{-1}{3}<\Gamma_3<1\), where they are saddles.  The sign of \(r(\Equilibrium_3)\) shows that \(\Equilibrium_3\) is a saddle for \(\frac{-5}{3}<\Gamma_3<1\), inside the deltoid, and a center outside it.

\begin{figure}
    \centering
    \includegraphics[width=\linewidth]{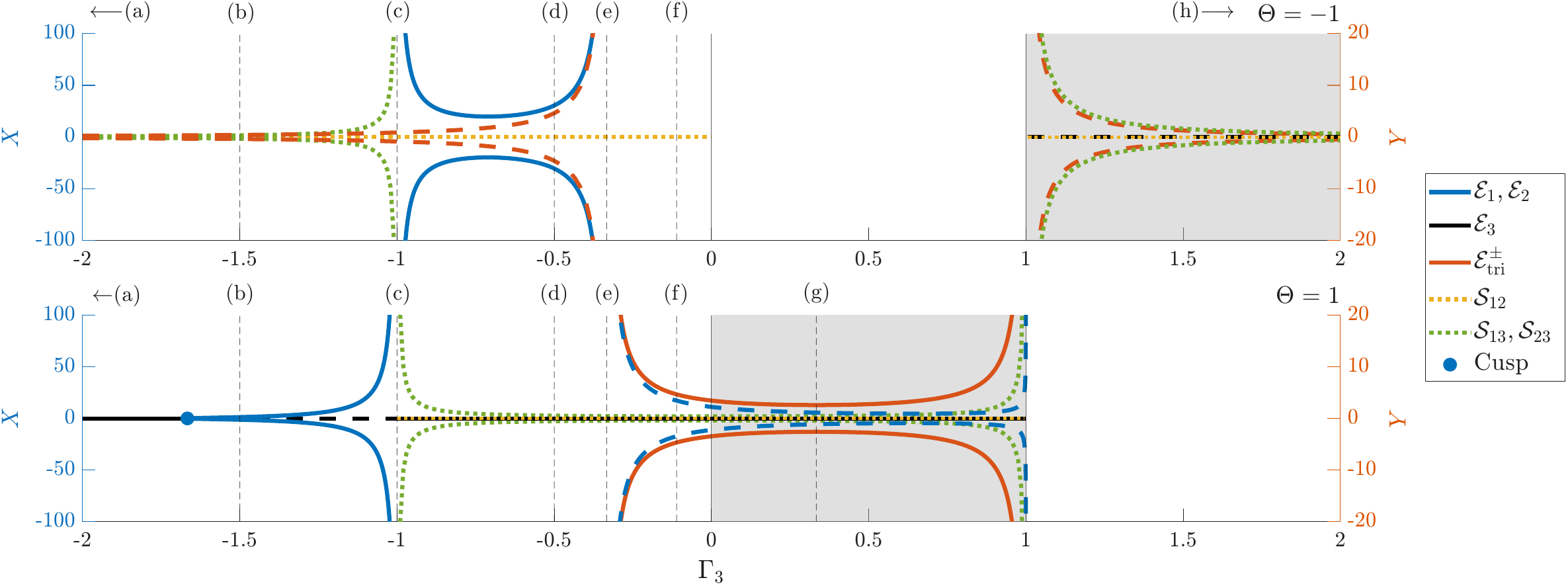}
    \caption{Bifurcation diagram displaying the equilibria and singularities defined by Eqs.~\eqref{XZ_singularities},~\eqref{equilateral_configs} and~\eqref{symmetric_equilibria} for \(\Theta=-1\) and \(\Theta=1\) as a function of \(\Gamma_3\) along the vertical center line of Fig.~\ref{fig:trilinear}. The \(X\)-component of the collinear equilibria and singularities is shown, using the \(y\)-axis scale on the left, and the \(Y\) component of the triangular equilibria is shown using the \(y\)-axis scale on the right. Dashed vertical lines show the \(\Gamma_1\) values indicated by the letter labels in Fig.~\ref{fig:trilinear}. Dashed equilibrium curves indicate instability. The points \textbf{a} and \textbf{h} lie outside the region plotted in this figure.}
    \label{fig:bifurcation_in_Gamma3}
\end{figure}

Fig.~\ref{fig:bifurcation_in_Gamma3} displays a bifurcation diagram for \(\Gamma_1=\Gamma_2\) and variable \(\Gamma_3\) with separate plots for the \(\Theta=-1\) and \(\Theta=1\) cases. The collinear and equilateral states are depicted, along with the singularities. The bifurcation at \(\Gamma_3=-1\) is a pitchfork, but all other bifurcations occur when fixed points diverge to and return from infinity, sometimes moving between the \(\Theta=\pm1\) surfaces.

\section{Phase portraits}
\label{sec:phase_portraits}

The dynamics occur on a symplectic leaf that may be either an ellipsoid or a hyperboloid, depending on the triplet of circulation values. In the former case, all trajectories are bounded since the leaf is compact. In the latter, we find that unbounded trajectories are possible only when the circulations satisfy an additional algebraic condition. Understanding the phase-plane dynamics also requires understanding how heteroclinic orbits connecting saddle points divide the leaves into non-overlapping regions. 

In this section, we plot phase portraits at a sequence of points labeled \textbf{a}--\textbf{h} along the line \(\Gamma_1 = \Gamma_2 =\frac{1-\Gamma_3}{2}\), the vertical symmetry axis in Fig.~\ref{fig:trilinear}.  Along this axis, the phase space is a hyperboloid for \(\Gamma_3<0\) and an ellipsoid for \(\Gamma_3>0\). In the former case, the phase surface is unbounded, but all orbits are bounded except for values of the circulations and \(\Theta\) satisfying additional algebraic conditions. We note the few cases where unbounded orbits exist. The phase surface is bounded in the latter, ellipsoidal case, as are all orbits. In addition to marking the equilibria and singularities, we plot the separatrix orbits emerging from any hyperbolic equilibria, and families of periodic orbits organized by these separatrices.  

Understanding the geometric meaning of the variables \(X\) and \(Y\) defined by Eq.~\eqref{muXYZ} can help us interpret these phase-space diagrams. First, consider Eq.~\eqref{jacobi_coordinates}. \(Z_1\) and \(Z_2\) may represent complex variables or vectors in \(\Reals^2\). We choose the latter.
Referring to Figure~\ref{fig:Jacobi}(a),  \(Z_1\) is the vector from \(z_2\) to \(z_1\) and \(\tilde{z}_2\) is their center of vorticity. The vector \(Z_2\) points from \(z_3\) to \(\tilde{z}_2\).  Then
\[
X = Z_1 \cdot Z_2 \qand Y = (Z_1 \times Z_2) \cdot \hat{k},
\]
where \(\hat{k}\) is the unit vector out of the plane. 
Thus, \(Y=0\)  when the vortices are collinear. When \(\Gamma_1 = \Gamma_2\), \(\tilde{z}_2\) lies at the midpoint of \(z_1\) and \(z_2\), so that \(X=0\) when the triangle is isosceles with the third vortex adjacent to the two equal sides. This last conclusion is only valid when \(\Gamma_1 = \Gamma_2\).

\subsection{Hyperboloidal Symplectic Leaves}

For \(\Gamma_3 < 0 \), the circulation parameters lie in the unshaded region of the diagram in Fig.~\ref{fig:trilinear} corresponding to hyperboloidal symplectic leaves, which degenerate to cones when \(\Theta=0\). We project this surface onto \(Z=0\) and plot the resulting trajectories as phase planes. All phase planes with \(\Theta<0\) are equivalent up to scaling. Those with \(\Theta>0\) form a second equivalence class. Equations~\eqref{XZ_singularities},~\eqref{equilateral_configs}, and~\eqref{symmetric_equilibria} by themselves are insufficient to determine which sheet a given singularity or equilibrium lies on. They must also be consistent with their signs, determined from~\eqref{ThetaMu} and~\eqref{muXYZ}.  Therefore, we plot three relevant phase planes, for~\(\Theta=-1\), \(0\), and~\(1\). 

The phase plane for \(\Theta=0\) is distinct, as the equilibria and singularities all degenerate to the triple collision \(\SingTriple\). In most cases, the topology of the trajectories is equivalent to the far-field topologies in the \(\Theta\neq0\) cases.  At the values plotted in Figs.~\ref{fig:point_c} and~\ref{fig:point_e}, the topology for \(\Theta=0\) does not follow this pattern.

\subsubsection*{Point \textbf{a}: \(\Gamma_3= \frac{-17}{3}\)}

\begin{figure}[ht]    \centering
   \includegraphics[width=\linewidth]{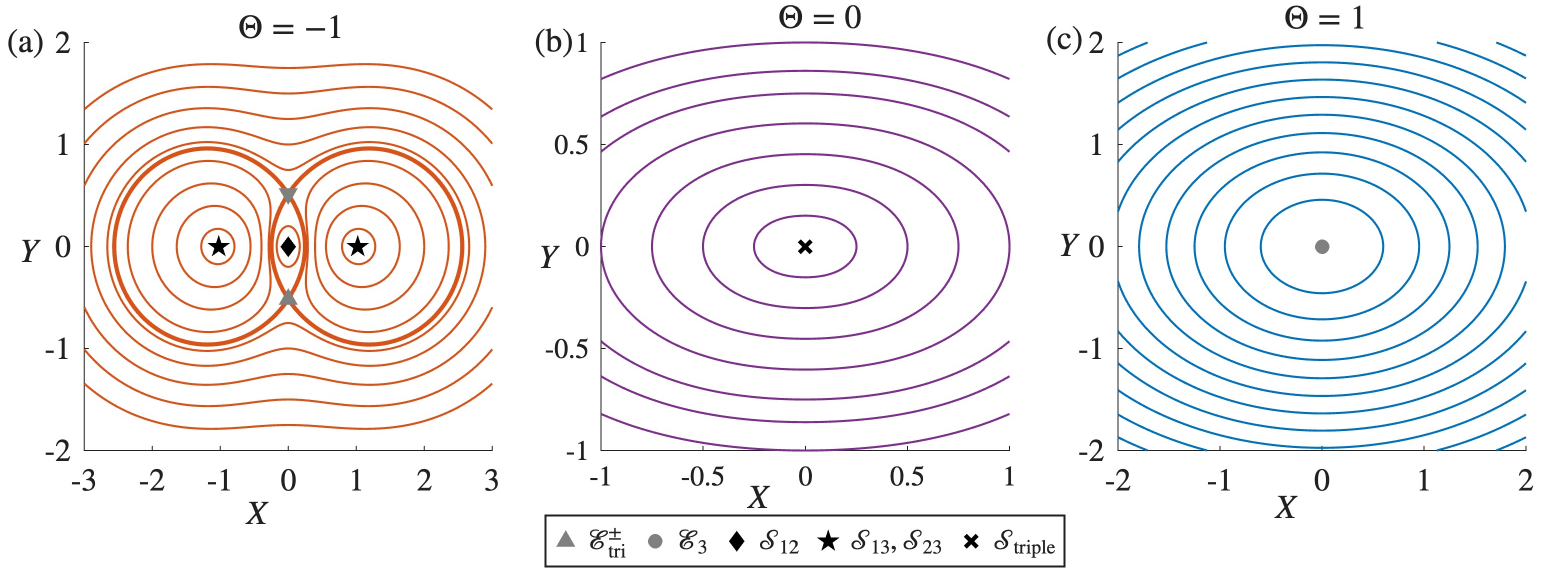}
   \caption{The \(XY\) phase planes at point \textbf{a} with \(\Gamma_3 =\frac{-17}{3}\). Thicker curves than the periodic orbits denote the separatrix orbits. The figures depict periodic orbits that are roughly equally spaced, rather than equally spacing the level sets of the Hamiltonian, which would lead to an accumulation of curves near each singularity. These conventions are used in all remaining figures.}
   \label{fig:point_a}
\end{figure}

Fig.~\ref{fig:point_a} shows phase planes or \(\Gamma_3< \frac{-5}{3}\). The leaves with \(\Theta<0\) contain all three singularities, and the two equilateral configurations, which are saddles. The leaves with \(\Theta>0\) contain a single equilibrium \(\Equilibrium_3\) in which the negatively-signed vortex 3 sits midway between the two positively-signed vortices. The \(\Theta=0\) leaf is topologically equivalent to the \(\Theta=1\) leaf with the singularity \(\SingTriple\) at its center.

\subsubsection*{Point \textbf{b}: \(\Gamma_3= \frac{-3}{2}\)}

\begin{figure}[ht]
    \centering
    \includegraphics[width=\linewidth]{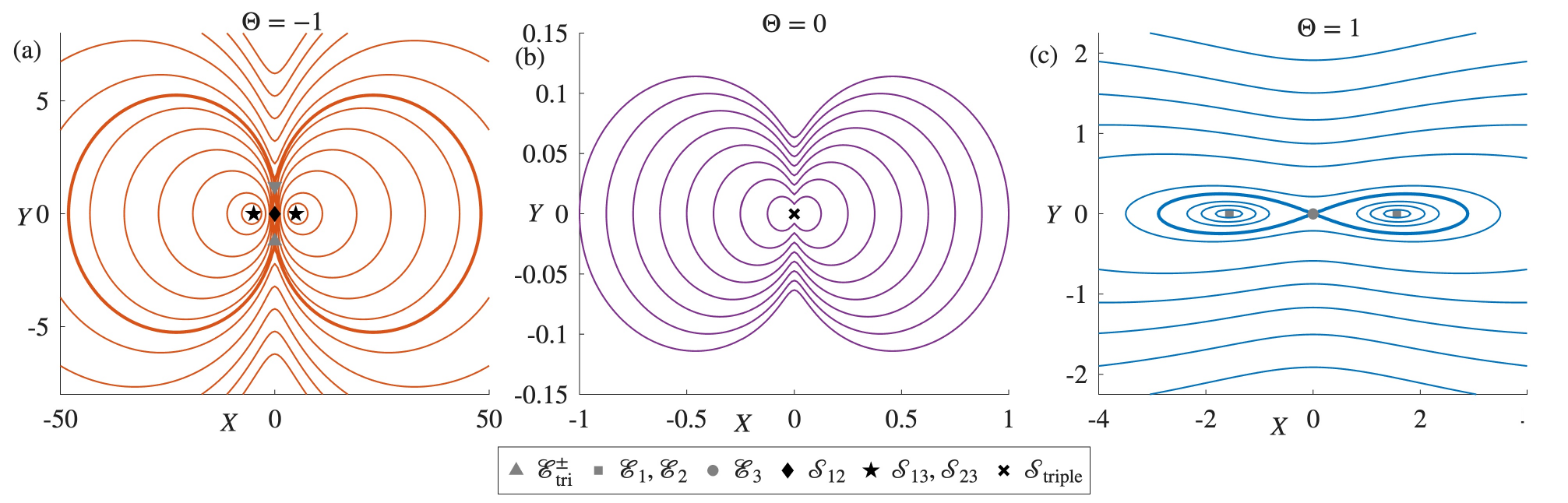}
    \caption{The \(XY\) phase planes of system at point \textbf{b} with \(\Gamma_3= \frac{-3}{2} \).}
    \label{fig:point_b}
\end{figure}

The phase surface for point \textbf{b} with \(\Gamma_3=\frac{-3}{2}\), shown in Fig.~\ref{fig:point_b}. The topology for \(\Theta<0\) is unchanged. At the same time, the collinear equilibria on the \(\Theta>0\) leaves have undergone a pitchfork bifurcation, with \(\Equilibrium_3\) becoming unstable, and the collinear equilibria \(\Equilibrium_1\) and \(\Equilibrium_2\) emerging as stable fixed points. The \(\Theta=0\) leaf is foliated by closed orbits.

\subsubsection*{Point \textbf{c}: \(\Gamma_3= -1\)}

As \(\Gamma_3 \nearrow -1\), the collinear equilibria \(\Equilibrium_1\) and \(\Equilibrium_2\) diverge to \(X=\pm \infty\), as do the singularities \(\Singularity_{13}\) and \(\Singularity_{23}\). Precisely at \(\Gamma_3=-1\), they all cease to exist. This is because at \(\Gamma_3=-1\), the third vortex may form a dipole with either of the other two vortices. Since the total circulation of the two vortices forming the dipole vanishes, their joint center of vorticity sits at infinity.

\begin{figure}[ht]
    \centering
    \includegraphics[width=\linewidth]{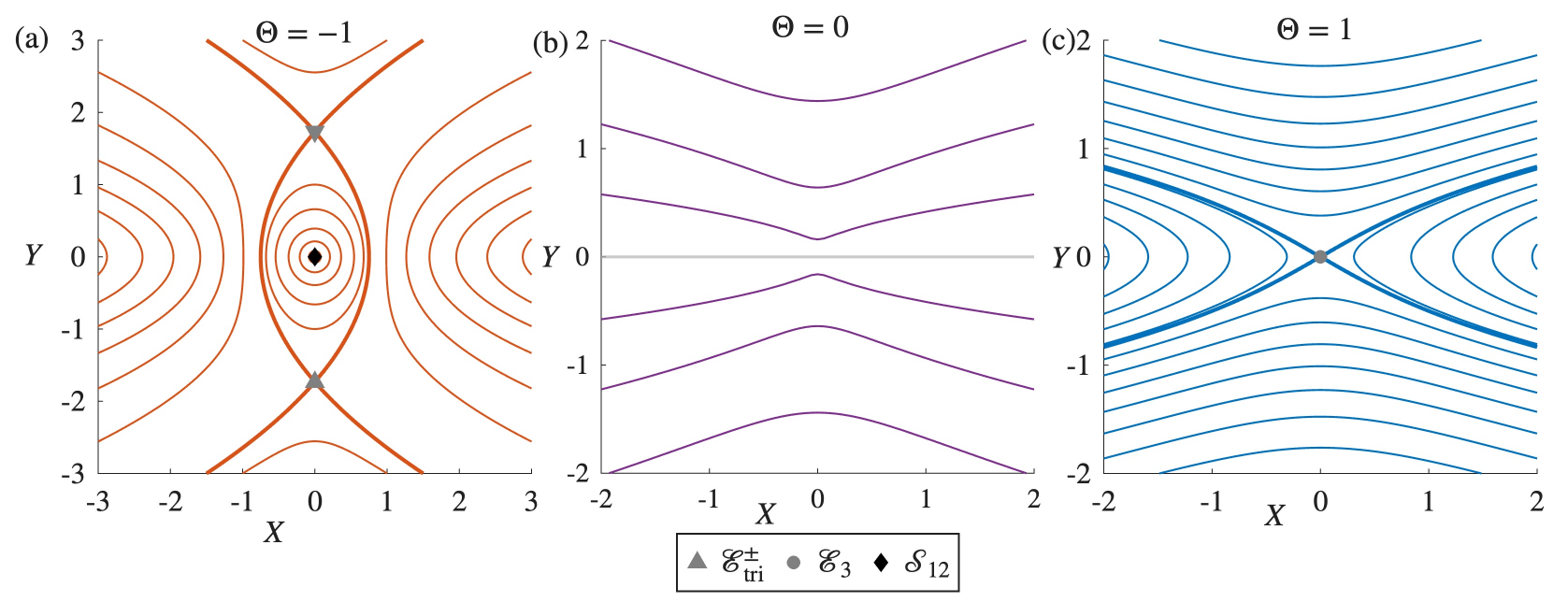}
    \caption{The \(XY\) phase planes at point \textbf{c} with \(\Gamma_3 =-1\). (a) The case \(\Theta<0\) with singularity \(\Singularity_{12}\) at the origin and triangular configurations \(\qty(X, Y, Z) \) at the intersections of the thick curves. (b) The case \(\Theta=0\). The gray line \(Y=0\) is singular. (c) The case \(\Theta > 0\) with collinear equilibrium at the origin.}
    \label{fig:point_c}
\end{figure}

This case features unbounded orbits on which two vortices form a dipole that scatters off a third, initially stationary vortex. In~\cite{Anurag:2024}, we used the phase surface plots in Fig.~\ref{fig:point_c} to derive a simplified explanation of this scattering phenomenon. In that paper, we analyze the scattering dynamics when the third vortex has arbitrary circulation. This occurs for circulations along either of the dashed lines through point \textbf{c} in Fig.~\ref{fig:trilinear}.

For this value of \(\Gamma_3\), when \(\Theta=0\), the three vortices must form a right triangle, a fact observed by Gröbli~\cite{Grobli:1877}. As such, they cannot be collinear, so \(Y\neq0\) and the \(X\)-axis is singular, as represented by the gray line in the center panel of Fig.~\ref{fig:point_c}.

\subsubsection*{Point \textbf{d}: \(\Gamma_3= \frac{-1}{2}\)}

\begin{figure}[ht]
    \centering
    \includegraphics[width=\linewidth]{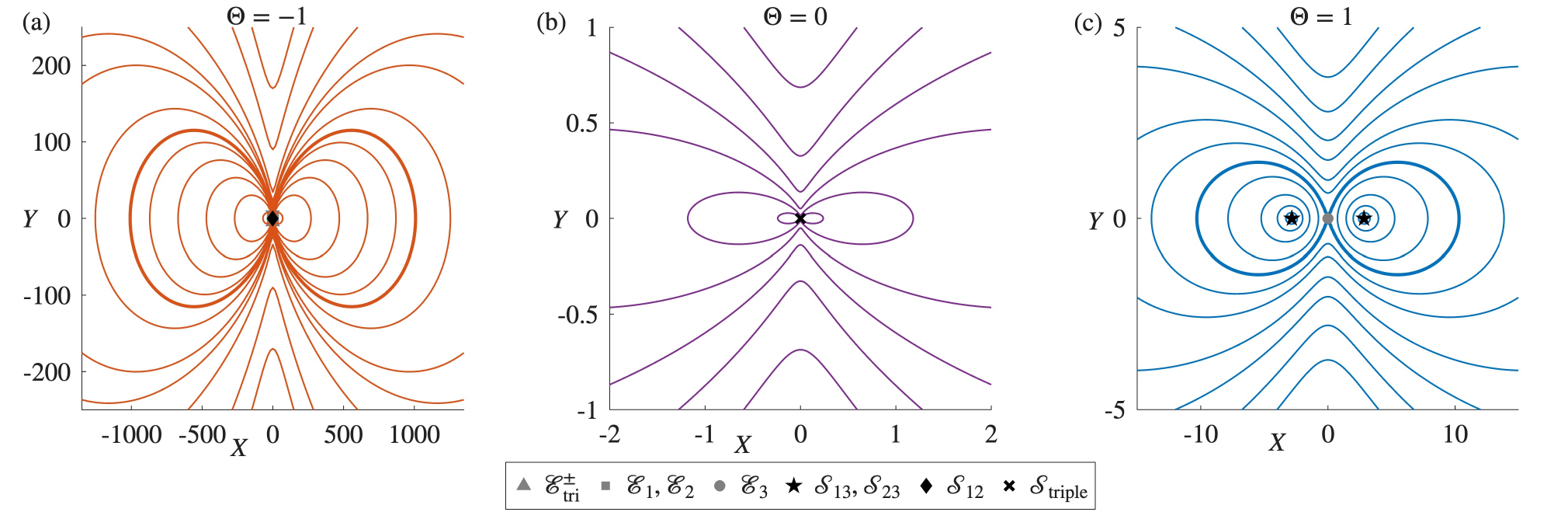}
    \caption{The \(XY\) phase planes of system at point \textbf{d}, \(\Gamma_3= \frac{-1}{2} \). The two equilibria \(\Equilibrium_1\) and \(\Equilibrium_2\) have moved to the \(\Theta<0\) surfaces and become centers. The two equilibria \(\EqTriPM\) have returned to the \(\Theta<0\) surfaces but have become centers. A close-up of the left image is shown in Fig.~\ref{fig:point_d_closeup}.}
    \label{fig:point_d}
\end{figure}

\begin{figure}[ht]
    \centering
    \includegraphics[width=0.4\linewidth]{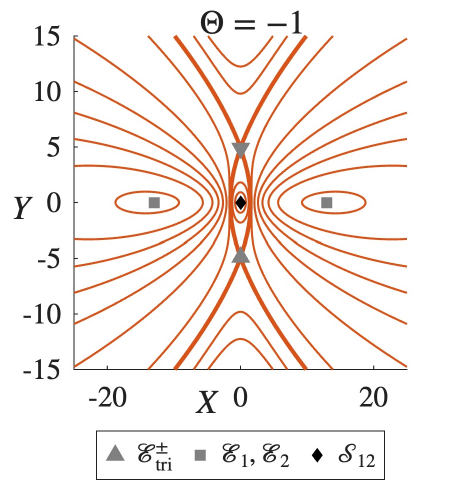}
    \caption{A closeup of the \(\Theta<0\) phase surface at point \textbf{d}, showing the singularities and equilibria. }
    \label{fig:point_d_closeup}
\end{figure}

\begin{figure}[ht]
    \centering
    \includegraphics[width=\linewidth]{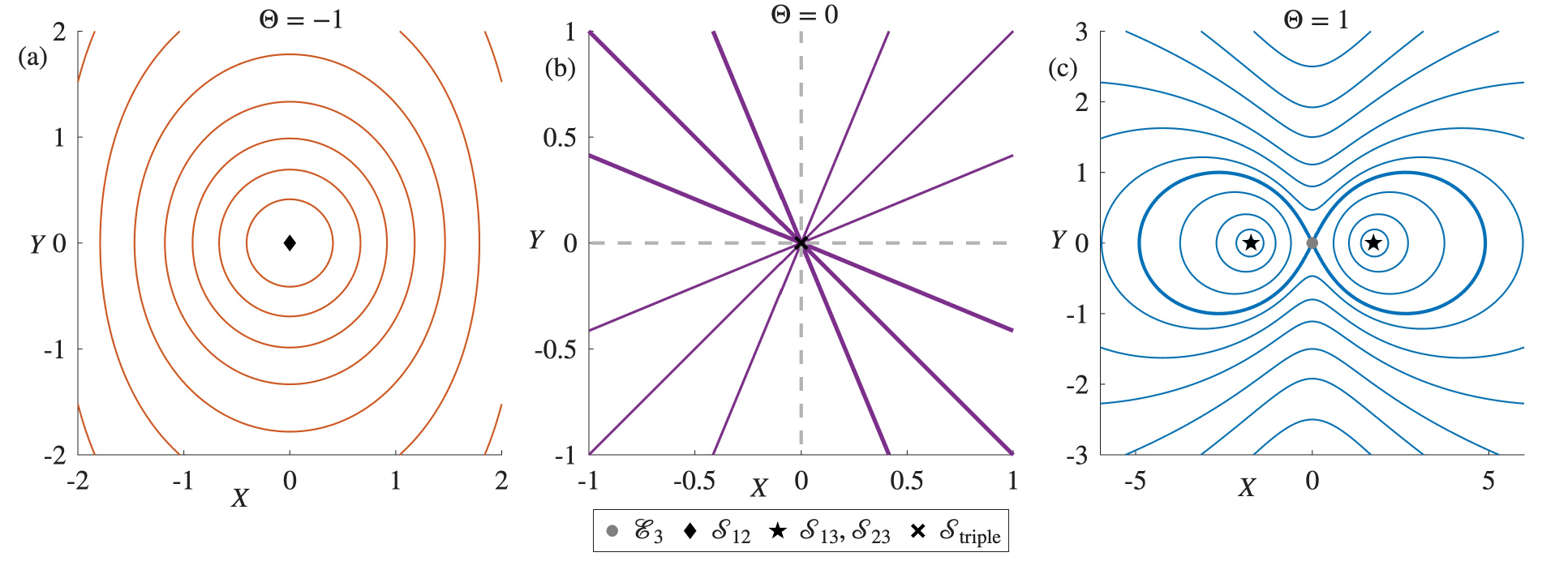}
    \caption{The \(XY\) phase planes at point \textbf{e} with \(\Gamma_3 =\frac{-1}{3}\). The trajectories in the \(\Theta=0\) case are rays through the origin and correspond to triangles that shrink to zero in finite forward time (quadrants one and three), finite backward time (quadrants two and four), or do not change size (the \(X\)- and \(Y\)-axes).}
  \label{fig:point_e}
\end{figure}

As \(\Gamma_3\) increases through \(-1\), the equilibria \(\Equilibrium_1\) and \(\Equilibrium_2\) return from infinity and reappear on the  \(\Theta<0\) leaves, remaining saddles, as shown in Fig.~\ref{fig:point_d}. At the same time, the singularities \(\Singularity_{13}\) and \(\Singularity_{23}\) return from infinity on the \(\Theta>0\) leaves.  Two pairs of heteroclinic orbits connect the saddles: one pair surrounds the singularity \(\Singularity_3\) at the origin for \(\Theta<0\), and a second, longer pair makes wide excursions enclosing the collinear equilibria \(\Equilibrium_1\) or \(\Equilibrium_2\). All orbits on the surface \(\Theta=0\) are closed. The \(\Theta>0\) phase surface is unchanged from the previous plot. Two heteroclinic orbits on the \(\Theta<0\) surfaces take such long excursions that the equilibria and singularities on this surface are not visible in the plot. Fig.~\ref{fig:point_d_closeup} shows a closeup highlighting this region.

\subsubsection*{Point \textbf{e}: \(\Gamma_3= \frac{-1}{3}\)}

Point \textbf{e} lies on the circle \(\gamma_2=0\) in Fig.~\ref{fig:trilinear}. As this point is crossed, the radius of the triangular equilibria diverges by Eq.~\eqref{equilateral_configs}, and Eq.~\eqref{charPolyTri} implies a change in their stability. The two collinear equilibria \(\Equilibrium_1\) and \(\Equilibrium_2\) also diverge to infinity, and their stability also changes. Therefore, exactly on this circle, the number of fixed points decreases from five to one.
Earlier literature has noted that \(\gamma_2=0\) is a necessary condition for finite-time collapse~\cite{Aref:1979, Grobli:1877, Krishnamurthy:2018}. This occurs when \(\Theta=0\), where the system has unbounded trajectories. The Hamiltonian takes the form
\begin{equation} \label{H_collapse}
    H = \frac{1}{9} \log{\left(\frac{4 Z^2-3 X^2}{Z^2}\right)}.
\end{equation}
Using Eq.~\eqref{rankoneXYZ} with \(\Theta=0\) to eliminate \(Z\) in favor of \(Y\) yields
\[
H = \frac{1}{9} \log{\left(\frac{ X^2+4 Y^2 }{X^2+Y^2}\right)}=\frac{1}{9}\log{\left( 5-3 \cos{2\theta}\right)}
\]
by the standard polar coordinate substitution. The level sets of \(H\) are rays through the origin since \(H\) is \(r\)-independent. The evolution reduces to
\[
\dv{r}{t}= \frac{2 \sqrt{3}\sin{2\theta}}{3\cos{2\theta}-5}, \, \dv{\theta}{t}=0.
\]
Thus, for \((X,Y)\) in the first or third quadrant, \(\dv{r}{t}\) is negative and constant. The solution shrinks to zero in finite time while maintaining a constant triangular profile. Similarly, in the second and fourth quadrants, the solution collapses to the origin at a finite negative time. The picture is essentially unchanged when \(\Gamma_2 \neq \Gamma_1\), except that the angles of the rays separating growing from decaying motions change.

The \(\Theta<0\) phase surface is foliated by closed orbits surrounding the singularity \(\Singularity_3\). The \(\Theta>0\) has one equilibrium \(\Equilibrium_3\), from which emerge a pair of homoclinic orbits, each encircling singularity \(\Singularity_1\) or \(\Singularity_2\).

\subsubsection*{Point \textbf{f}: \(\Gamma_3= \frac{-1}{9}\)}

\begin{figure}[htp]
    \centering
    \includegraphics[width=\linewidth]{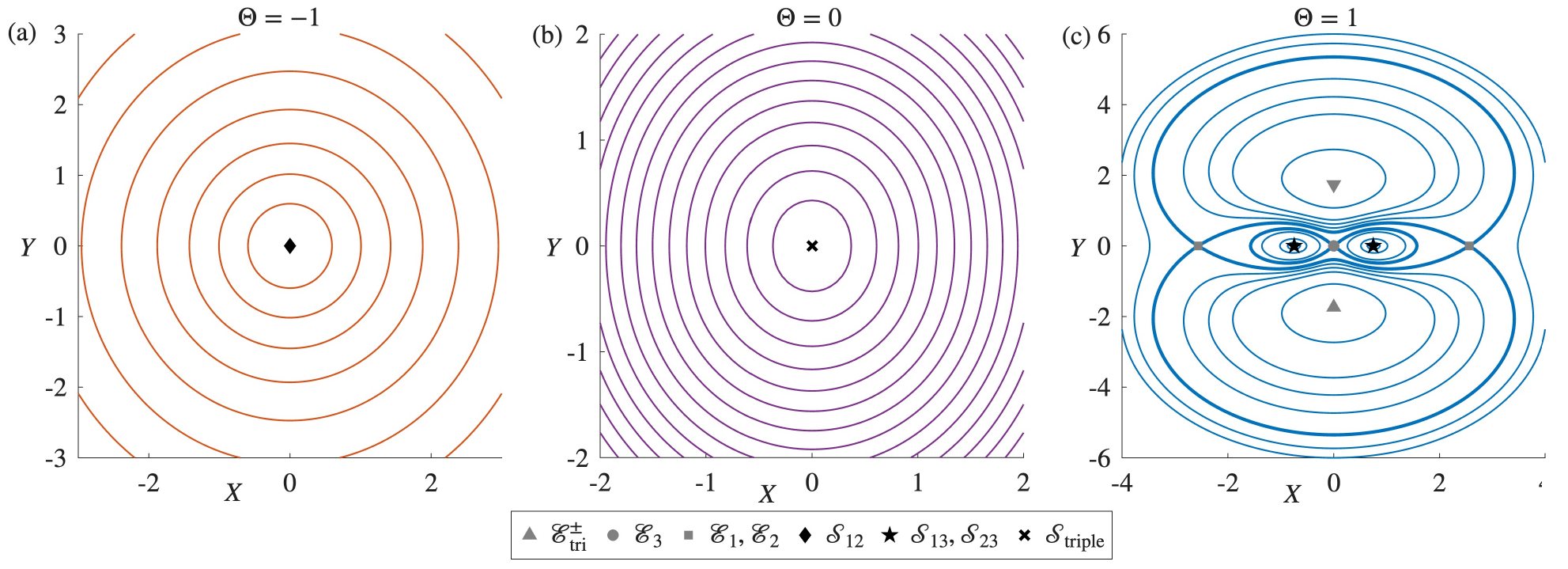}
    \caption{The \(XY\) phase planes in case~\textbf{f} with \(\Gamma_3=\frac{-1}{9}\). (a) The case \(\Theta<0\). (b) The case \(\Theta=0\). (c) The case \(\Theta>0\). }
 \label{fig:point_f}
\end{figure}

As \(\Gamma_3\) increases through \(-\frac{1}{3}\), the equilibria \(\EqTriPM\) diverge to infinity, and re-emerge on the leaves with \(\Theta>0\), now as centers. Simultaneously, the equilibria \(\Equilibrium_1\) and \(\Equilibrium_2\) diverge to infinity and return on the \(\Theta>0\) leaves, now as saddle points. 
The phase plane at point \textbf{f} with \(\Gamma_3 = \frac{-1}{9}\) is shown in Fig.~\ref{fig:point_f}.  It features the two stable equilateral equilibria, the three unstable collinear equilibria, and the singularities \(\Singularity_{13}\) and \(\Singularity_{23}\), and six heteroclinic orbits connecting the collinear equilibria. The dynamics for \(\Theta<0\) is foliated by periodic orbits surrounding the singularity \(\Singularity_{12}\).

\subsection{Ellipsoidal Symplectic Leaves}

As \(\Gamma_3\) increases through 0, the quadric form defining the symplectic leaves changes signature, leading to an ellipsoidal phase space.  We plot the ellipsoidal symplectic leaves as spheres, ignoring that their axes are unequal, plotting \(Y=0\) as the equator and the north and south poles at \(Y=\pm 1\).

  \subsubsection*{Point \textbf{g}: \(\Gamma_3 =\frac{1}{3}\)}

\begin{figure}[htb]
 \centering
   \includegraphics[width=0.8\linewidth]{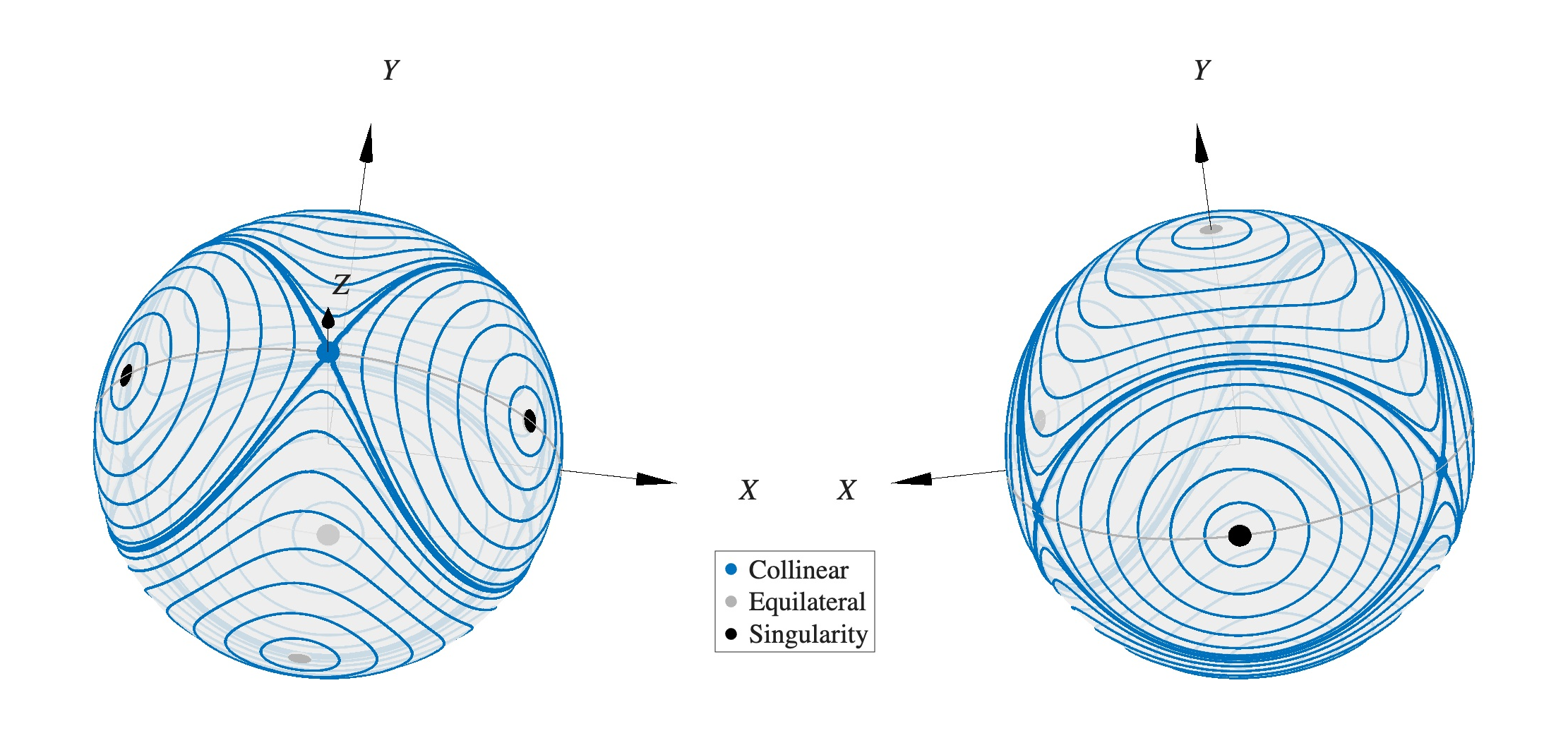}
   \caption{The phase portrait at point \textbf{h} with \(\Gamma_3 =\frac{1}{3}\), showing "front" and "back" views of the sphere.}
     \label{fig:sphere_h}
 \end{figure}

Other than the topology, the only change in the phase space of cases \textbf{f} and \textbf{g} is that the singularity \(\Singularity_{12}\) has moved to the \(\Theta>0\) leaves. Setting \(\Gamma_3=\frac{1}{3}\) yields the most symmetric case where all three vortices have identical circulation, the case represented by the point \textbf{h} in Fig.~\ref{fig:trilinear}.  Its phase portrait is shown in Fig.~\ref{fig:sphere_h}. The three collinear equilibria lie equispaced on the equator. Halfway between each pair lies a singularity; for example \(\Singularity_{12}\) is between equilibria \(\Equilibrium_1\) and \(\Equilibrium_1\). The equilateral equilibria \(\EqTriPM\) sit at the north and south poles

The six heteroclinic orbits divide the sphere into five families of periodic orbits. Each heterocline connects two equilibria, corresponding to an orbit that exchanges a vortex from the end of the line segment with the one at its center. A family of periodic orbits surrounds each of the singularities on the equator. The limiting orbits of these families are heteroclinic cycles connecting pairs of collinear relative equilibria. An additional family of periodic orbits surrounds each of the equilateral equilibria. The limiting orbits of these families are heteroclinic chains that visit all three collinear equilibria.

 \subsubsection*{Point \textbf{g}: \(\Gamma_3 > 1\)}
\begin{figure}[htbp]
     \centering
      \centering
   \includegraphics[width=0.8\linewidth]{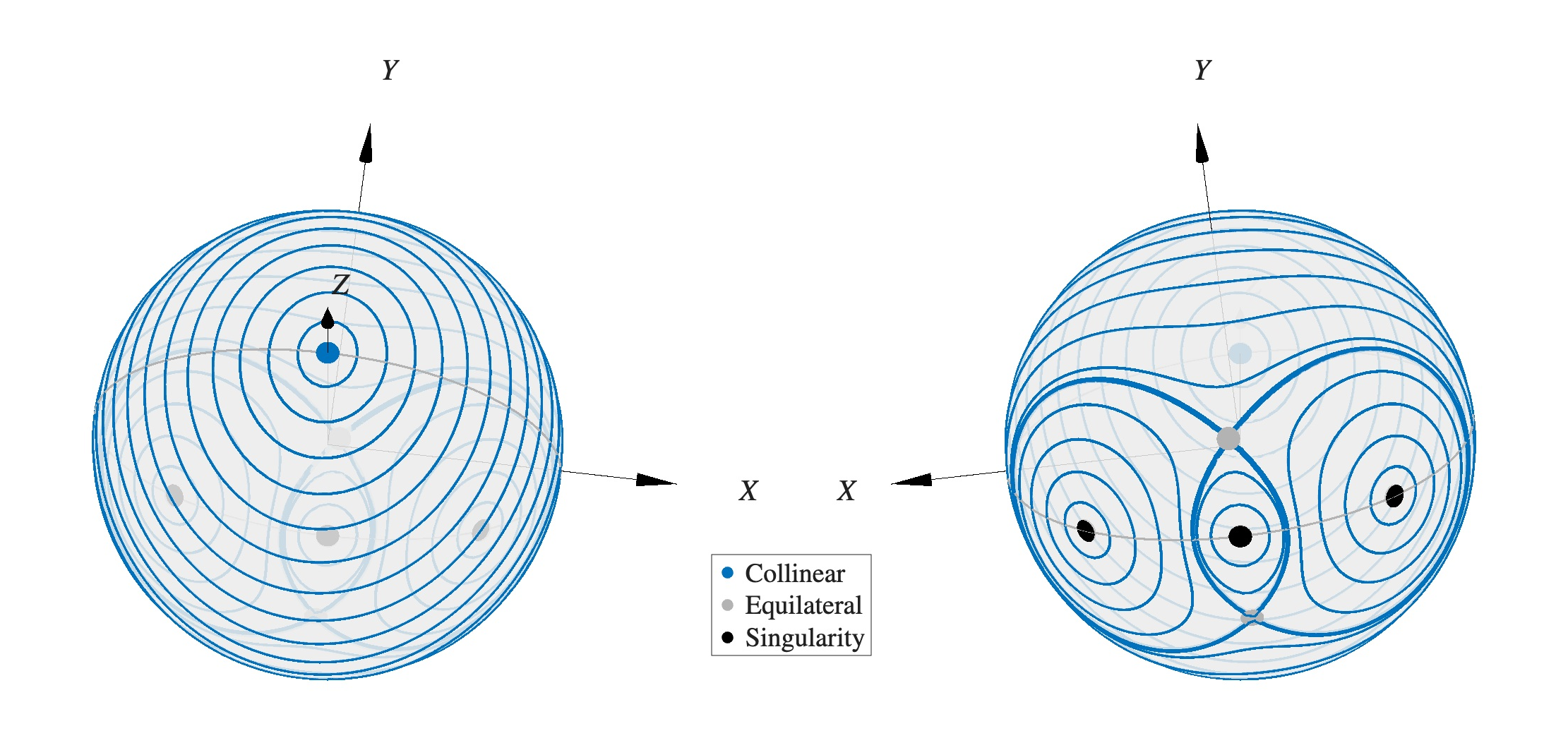}
     \caption{The phase portrait at point \textbf{h} with \(\Gamma_3=5\). Only the single collinear equilibrium \(\Equilibrium_3\) remains and is a center. The two equilateral triangular equilibria \(\EqTriPM\) are saddles.}
     \label{fig:sphere_j}
 \end{figure}

Finally, The phase portraits for \(\Gamma_3 = 5\) is shown in Fig.~\ref{fig:sphere_j}. At \(\Gamma_3=1\), \(\Gamma_1=\Gamma_2=0\), so the motion degenerates to that of a single vortex. As \(\Gamma_3\) crosses 1, the parameters move from the central shaded triangle in Fig.~\ref{fig:trilinear} to the unbounded shaded region above it. The collinear equilibria \(\Equilibrium_1\) and \(\Equilibrium_2\) disappear, since they only exist inside the deltoid region. The triangular equilibria \(\EqTriPM\) change from centers to saddles, as this bifurcation occurs on the circle. Heteroclinic cycles connecting \(\EqTriPM\) encircle each of the three collinear singularities and the families of periodic orbits surrounding them. The stability of the remaining collinear equilibrium has also changed from unstable to stable, and it is surrounded by periodic orbits.

\section{Conclusion}
\label{sec:conclusion}
We have combined a Jacobi coordinate reduction, a Lie-Poisson reduction, following Ohsawa, and a judicious choice of coordinates to derive a reduced set of evolution equations for the motion of three vortices in the plane. The combined reduction introduces no singularities into the coordinate system, in contrast with formulations based on the reduced system Gröbli used to integrate the system.

Using this reduction, we studied the global phase portraits of the three-vortex problem, the bifurcations of relative equilibria, and their stability. In particular, collinear configurations, which are singular points in the barycentric formulation, are regular points in the current reduction. This allows us to analyze their linear stability using the standard linearization method.

This method of analysis should be useful in other contexts. 
The barycentric formulation has also been used to study two other integrable cases of vortex motion, namely the motion of three vortices on a sphere and the motion of four vortices with vanishing total circulation and linear impulse~\cite{Kidambi.1998, Aref:1999, Tsai:2020}. The reduction method described here will allow us to map out the dynamics of these systems more comprehensively than has been done to date.

\appendix

\section{Vanishing total circulation}
\label{sec:vanishing}

The last step of the Jacobi reduction cannot be applied when the total circulation vanishes, so we must use an alternative procedure. We may always label the vortices so that the sequence of transformations works at each step except the last. At the final step, we need an alternative to the reduction described by Eqs.~\eqref{Jacobi1} that applies in the dipole case when \(\Gamma_1+\Gamma_2=0\). Ohsawa has noted the need for a different reduction method when the total circulation vanishes~\cite{Ohsawa:2024}. This can be traced to the \emph{superintegrability} of the three-vortex problem with vanishing total circulation~\cite{Galajinsky.2022}.

\subsection{The reduction method}
Consulting Fig.~\ref{fig:Jacobi}(b), we define two position coordinates \(Q_j\) and momentum vectors \(P_j\)
\begin{equation} \label{jacobiZero}
\begin{aligned}
Q_1 &= \frac{x_1 + x_2}{2},  &Q_2 &= \frac{y_1 + y_2}{2}, \\
P_1  & = y_1 - y_2, & P_2 & = - x_1  + x_2.
\end{aligned}
\end{equation}
The Poisson bracket, and thus the evolution equations, maintain the dependence on the circulation:
\[
\dv{Q_j}{t} = \frac{1}{\Gamma_j}\pdv{H}{P_j}, \,
\dv{P_j}{t} = \frac{-1}{\Gamma_j}\pdv{H}{Q_j}.
\]
For the case \(N=2\), the Hamiltonian becomes 
\[
H = \frac{\Gamma_1^2}{2}\log{(P_1^2 + P_2^2)}.
\]

The momenta are conserved since \(H\) is \(Q_j\)-independent, and the particles move in a straight line at constant velocity 
\[
\dv{t} \binom{Q_1}{Q_2}=\frac{\Gamma_1}{P_1^2+P_2^2}\binom{P_1}{P_2}.
\]
This is somewhat different than the Hamiltonian formulation~\eqref{N_vortex_ham_eqns}, in which the \(x\) and \(y\) coordinates play the roles of position and momentum. Here, we have returned to the more standard case in which \(Q\) variables represent positions and \(P\) variables momenta. 

\begin{figure}[ht]
    \centering
    \includegraphics[width=0.33\linewidth]{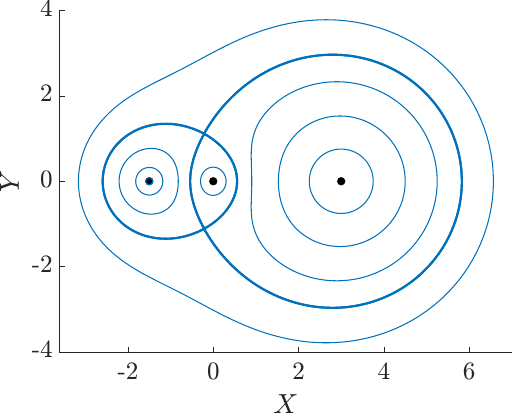}
    \caption{Thehase plane with zero total circulation, with \((\Gamma_1,\Gamma_2,\Gamma_3)=\left(2,1,-3\right)\).}
    \label{fig:zero_circulation}
\end{figure}

We now consider the case of three vortices with \(\gamma_1=0\), and, more particularly, \(\Gamma_3 = -\Gamma_1 - \Gamma_2 \neq 0\). We may apply the reduction~\eqref{Jacobi1} to vortices 1 and 2, but the second recursive application of the Jacobi coordinate reduction cannot be used to combine the virtual vortex at position \(\tilde{z}_2\) with vortex 3. Instead, we apply the transformation~\eqref{jacobiZero} to these coordinates. Without loss of generality, we may take \(P_2=0\). If \(P_3 \neq 0\), we may scale it to one without loss of generality. Letting \(Q_1=X\) and \(P_1=Y\) and ignoring additive constants, this yields a Hamiltonian
\begin{equation} \label{HXY}
\begin{split}
H=& -\frac{\Gamma_1 \Gamma_2}{2}  \log \left(X^2+Y^2\right)
+\frac{ \left(\Gamma_1+\Gamma_2\right)\Gamma_1}{2}  \log \left(\left(\Gamma_1+\Gamma_2 (X-1)\right)^2+\Gamma_2^2 Y^2\right)\\
&+\frac{\left(\Gamma_1+\Gamma_2\right) \Gamma_2}{2}  \log \left(\left(\Gamma_2+\Gamma_1 (X+1)\right)^2+\Gamma_1^2 Y^2\right).
\end{split}
\end{equation}
The three singularities are at \(Y=0\) and 
\[
X_{12} = 0, \quad
X_{13} = 1-\frac{\Gamma_1}{\Gamma_2}, \quad 
X_{23} = -1-\frac{\Gamma_2}{\Gamma_1}.
\]
Rott and Aref derived this system separately by different methods, as did Behring~\cite{Aref:1989, Rott:1989, Behring:2020}. It has equilateral relative equilibria at 
\begin{equation} \label{EqTriPM0}
   \EqTriPM  = \parentheses
{\frac{\Gamma_2^2-\Gamma_1^2}{2 \left(\Gamma_1^2+\Gamma_2 \Gamma_1+\Gamma_2^2\right)},
\pm\frac{\sqrt{3} \left(\Gamma_1+\Gamma_2\right)^2}{2 \left(\Gamma_1^2+\Gamma_2 \Gamma_1+\Gamma_2^2\right)}},
\end{equation}
A standard computation shows that these are saddles, connected by heteroclinic orbits. One such phase plane is shown in Fig.~\ref{fig:zero_circulation}. 

When \(P_2=P_3=0\), the system is super-integrable, and the Hamiltonian reduces further to 
\begin{equation} \label{Hnothing}
H = \frac{\Gamma_1^2+\Gamma_2 \Gamma_1+\Gamma_2^2}{2}  \log \left(X^2+Y^2\right),
\end{equation}
again, ignoring an additive constant. The singularity at the origin corresponds to a triple collision.

\section{The resultant and the discriminant}
\label{sec:resultant}

Aref gives an excellent introduction to the resultant and discriminant and applies it to this problem in Ref.~\cite{Aref:2009}. Mathematica has functions to compute both.

\subsubsection*{The resultant}
The resultant must vanish for two polynomials \(a(x)\) and \(b(x)\) to have a common root. Consider, for example, the case that \(a(x) = a_2 x^2 + a_1 x + a_0\) and \(b(x) = b_3 x^3 + b_2 x^2 + b_1 x + b_0\).
Augment the system 
\[
a(x)=0, \quad b(x) = 0
\]
with the redundant equations 
\[
x^2 a(x) = 0, \quad x a(x) = 0, \qand x b(x) = 0.
\]
These may be arranged into a linear system
\begin{equation} \label{sylvester}
\begin{bmatrix}
a_2 & a_1 & a_0 &   0 &   0 \\
  0 & a_2 & a_1 & a_0 &   0 \\
  0 &   0 & a_2 & a_1 & a_0 \\
b_3 & b_2 & b_1 & b_0 &   0 \\
  0 & b_3 & b_2 & b_1 & b_0
\end{bmatrix}
\begin{bmatrix}x^4 \\ x^3 \\ x^2 \\ x^1 \\ x^0 \end{bmatrix} 
=
\begin{bmatrix} 0 \\ 0 \\ 0 \\ 0 \\ 0 \end{bmatrix}.
\end{equation}
A necessary condition for system~\eqref{sylvester} to have a solution is the vanishing of its determinant, the \emph{Sylvester determinant}. This is the \emph{resultant} of the two polynomials \(\res{p}{q}{x}\). The above construction generalizes straightforwardly to polynomials of arbitrary degree. The dimension of the Sylvester matrix is \((\deg{a}+\deg{b})\), so the size and complexity of the resultant expression grow rapidly with the degrees of the polynomials. Several equivalent definitions exist, and computing the determinant by multiplying out the terms is not an efficient algorithm. 

The resultant can be used to find the roots of a system of two polynomials in two variables \(p(x,y)\) and \(q(x,y)\). Treating each as a polynomial in \(y\) with \(x\)-dependent coefficients, the resultant eliminates \(y\) and returns a single higher-order polynomial in \(x\), which we can treat using standard methods. The degree of this polynomial then gives the number of (possibly complex-valued) roots. 

\subsubsection*{The discriminant}
At values of \(\mu\) where the number of roots of a polynomial \(p(x,\mu)\) changes, i.e., at bifurcation points, \(p(x,\mu)\) and \(p'(x,\mu)\) vanish simultaneously. The \emph{discriminant} of \(p\) is proportional to the resultant of \(p\) and \(p'\). It generalizes the \(b^2-4ac\) term in the quadratic formula, whose vanishing indicates a double root of the quadratic equation.

\section{Barycentric coordinates for the plane}
\label{sec:trilinear}

We introduce barycentric coordinates to visualize the parameter space of circulations. Barycentric coordinates provide a method to parameterize \(\Reals^n\) based on a set of points \(\{\sfP_1,\ldots,\sfP_{n+1} \}\) that form an \(n\)-simplex. Here we consider the case of three non-collinear points \(\sfP_j=(\sfx_j,\sfy_j) \in \Reals^2\).

Any point \(\sfP \in \mathbb{R}^2\) can be written uniquely as
\begin{equation} \label{wP}
\sfP = w_1 \sfP_1 + w_2 \sfP_2 + w_3 \sfP_3,
\text{ where }
w_1 + w_2 + w_3 = 1.    
\end{equation}
The \emph{weights} \((w_1,w_2,w_3)\) are the barycentric coordinates of \(\sfP\) with respect to \((\sfP_1, \sfP_2, \sfP_3)\), an idea introduced by Möbius~\cite{Mobius1827} and widely used in geometry and geometric modeling; see, e.g.,~\cite{Farin.2007}.

Written in vector form, Eq.~\eqref{wP} is
\[
\begin{pmatrix}
\sfx_1 & \sfx_2 & \sfx_3 \\
\sfy_1 & \sfy_2 & \sfy_3 \\
1   & 1   & 1
\end{pmatrix}
\begin{pmatrix}
   w_1 \\ w_2 \\  w_3
\end{pmatrix}
=
\begin{pmatrix}
   \sfx \\ \sfy \\ 1 
\end{pmatrix}
,
\]
so uniqueness follows from the matrix being invertible, which holds since the three points are non-collinear.

Geometrically, these coordinates describe the position of \(\sfP\) relative to the triangle \(T\) with vertices \(\sfP_1\), \(\sfP_2\), and \(\sfP_3\). If all three weights are positive, then \(\sfP\) lies inside \(T\). If one coordinate is zero and the others are nonnegative, then \(\sfP \in \partial T\). If one or more coordinates are negative, then \(\sfP\) lies outside \(T\).

Here, the condition~\eqref{circulation_1} allows us to interpret the circulations as weights. Once we choose three points \(\sfP_j\), then each triplet of circulations corresponds to a unique point in the plane. Following~\cite{Conte:1979, Aref:2009}, we take three points equally spaced around the unit circle:
\[
\sfP_1 = \left(\frac{-\sqrt{3}}{2},\frac{-1}{2}\right), \,
\sfP_2 = \left(\frac{\sqrt{3}}{2},\frac{-1}{2}\right), \,
\sfP_3 = \left(0,1\right).
\]
Thus, the triple of vorticities gets mapped to the point
\begin{equation} \label{barycenter_explicit}
\parentheses{\sfx,\sfy} = \parentheses{
\frac{\sqrt{3}}{2}\left(-\Gamma_1+\Gamma_2\right),
-\frac{\Gamma_1}{2}-\frac{\Gamma_2}{2}+\Gamma_3}.
\end{equation}

Relatedly, Aref used barycentric coordinates to analyze Gröbli’s system~\eqref{l_eqn} as follows~\cite{Aref:1979}. The conservation laws~\eqref{constants_of_motion} imply that 
\[ 
\Gamma_1 \Gamma_2 \ell_{12}(t)^2 + 
\Gamma_2 \Gamma_3 \ell_{23}(t)^2 + 
\Gamma_3 \Gamma_1 \ell_{31}(t) ^2  
= \Gamma_1 \Gamma_2 \Gamma_3 C 
\]
defines a constant \(C\). Defining the rescaled squared lengths
\[
b_1=\frac{\ell_{23}^2}{\Gamma_1 C}, \quad b_2=\frac{\ell_{31}^2}{\Gamma_2 C}, \quad b_3=\frac{\ell_{12}^2}{\Gamma_3 C},
\]
then \(b_1+b_2+b_3=1\) This way, \(b_j\) and the solutions to system~\eqref{l_eqn} may be interpreted as barycentric coordinates for \(\Reals^2\). The requirement that the \(\ell_{ij}\) satisfy the triangle inequality leads to singularities in this representation whenever the three vortices are collinear, an impetus for the current project.

\section*{Acknowledgments}{This work was supported by the NSF under DMS-2206016. We thank Jared Bronski, Luis García-Naranjo, and Tomoki Ohsawa for useful discussions. Our translation of Gröbli's thesis depended on Google Translate and the Mathpix app~\cite{mathpix}. The paper benefited from the anonymous reviewers' incisive comments.}

\end{document}